\documentclass[11pt]{article}

\usepackage{graphicx}
\usepackage{subfigure}
\usepackage{color}
\usepackage{epsfig}
\usepackage{amssymb}
\usepackage{latexsym}

\newcommand {\C} {{\rm I\kern-5.5pt C}}

\newcommand{\bP}[1]{{\mathbb{P}}\left[{#1}\right]}
\newcommand{\bE}[1]{{\mathbb{E}}\left[{#1}\right]}

\newcommand{\1}[1]{{\bf 1}\left[#1\right]}       

\newcommand{\fsquare}{\vrule height6pt width7pt depth1pt}   
\newcommand{\myproof}{{\hfill \\ \bf Proof. \ }}           
\newcommand{\myendpf}{\hfill \fsquare \\[0.1in]}             

\newtheorem{theorem}{Theorem}[section]

\newtheorem{lemma}[theorem]{Lemma}
\newtheorem{proposition}[theorem]{Proposition}

\begin{document}

\title{A zero-one law 
       for the existence of triangles in random key graphs}

\author{
Osman Ya\u{g}an and Armand M. Makowski\footnote{This 
work was supported by NSF Grant CCF-07290.}
\\
{\tt oyagan@umd.edu}, {\tt armand@isr.umd.edu} \\
Department of Electrical and Computer Engineering\\
and Institute for Systems Research \\
University of Maryland, College Park, MD 20742.\\
}

\maketitle

\begin{abstract}
\normalsize 
Random key graphs are random graphs induced by the random key
predistribution scheme of Eschenauer and Gligor under the
assumption of full visibility. For this class of random
graphs we show the existence of a zero-one
law for the appearance of triangles,
and identify the corresponding critical scaling.
This is done by applying the method of first and second moments
to the number of triangles in the graph.
\end{abstract}


\section{Introduction}

Random key graphs are random graphs that belong to the class of random
intersection graphs \cite{SingerThesis}; in fact
they are sometimes called 
uniform random intersection graphs by some authors
\cite{GodehardtJaworski, GodehardtJaworskiRybarczyk}.
They have appeared
recently in application areas as diverse as clustering analysis
\cite{GodehardtJaworski, GodehardtJaworskiRybarczyk},
collaborative filtering in recommender systems \cite{Marbach2008}
and random key predistribution for wireless sensor networks (WSNs)
\cite{EschenauerGligor}. In this last context, random key graphs
naturally occur in the study of the following random key
predistribution scheme introduced by Eschenauer and Gligor
\cite{EschenauerGligor}: Before deployment, each sensor in a WSN
is independently assigned $K$ distinct cryptographic keys which
are selected at random from a pool of $P$ keys. These $K$ keys
constitute the key ring of the node and are inserted into its
memory. Two sensor nodes can then establish a secure link between
them if they are within transmission range of each other and if
their key rings have at least one key in common; see
\cite{EschenauerGligor} for implementation details. If we assume
{\em full visibility}, namely that nodes are all within
communication range of each other, then secure communication
between two nodes requires only that their key rings share at
least one key. 
The resulting notion of adjacency defines the
class of random key graphs; see Section \ref{sec:RandomKeyGraph}
for precise definitions.

Much efforts have recently been devoted to developing zero-one
laws for the property of connectivity in random key graphs. A key
motivation can be found in the need to obtain conditions under
which the scheme of Eschenauer and Gligor guarantees secure
connectivity with high probability in large networks. An
interesting feature of this work lies in  the following fact:
Although random key graphs are {\em not} equivalent to the
classical Erd\H{o}s-R\'enyi graphs \cite{ER1960}, it is possible
to transfer well-known zero-one laws for connectivity in
Erd\H{o}s-R\'enyi graphs to random key graphs by asymptotically
matching their edge probabilities. This approach, which was
initiated by Eschenauer and Gligor in their original analysis
\cite{EschenauerGligor}, has now been validated rigorously; see
the papers \cite{BlackburnGerke,
DiPietroManciniMeiPanconesiRadhakrishnan2008, Rybarczyk2009,
YaganMakowskiISIT2009, YaganMakowskiConnectivity} for recent
developments. Furthermore, Rybarczyk \cite{Rybarczyk2009} has
shown that this transfer from Erd\H{o}s-R\'enyi graphs also works
for a number of issues related to the giant component and its diameter.

In view of these successes, it is natural to wonder whether the
transfer technique can be applied to other graph properties. In
particular, in the literature on random graphs there is a long
standing interest \cite{ER1960, JansonLuczakRucinski,
KaronskiScheinermanSingerCohen, PenroseBook, SingerThesis} in the
containment of certain (small) subgraphs, the simplest one being
the {\em triangle}. This last case has some practical relevance
since the number of triangles in a graph is closely related to its
clustering properties \cite{YaganMakowskiChennai2009}.
With this in mind, we study the zero-one
law for the existence of triangles in random key graphs and
identify the corresponding critical scaling.

From these results we easily conclude
that in the many node regime,
the expected number of triangles
in random key graphs is always at least as
large as the corresponding quantity
in asymptotically matched Erd\H{o}s-R\'enyi graphs.
For the parameter range that is of
practical relevance in the context of WSNs,
this expected number of triangles
can be orders of magnitude larger
in random key graphs than in Erd\H{o}s-R\'enyi graphs,
a fact also observed earlier via simulations in
\cite{DiPietroManciniMeiPanconesiRadhakrishnan2008}. As a result,
transferring results from Erd\H{o}s-R\'enyi graphs by matching
their edge probabilities is not a valid approach in general, and
can be quite misleading in the context of WSNs.

The zero-one laws obtained here 
were announced in 
the conference paper \cite{YaganMakowskiAllerton2009}.
The results are established by making use of
the method of first and second moments to the
number of triangles in the graph.
As the discussion amply shows, the technical details, especially
for the one-law, are quite involved, and an outline 
of the proofs can be found in \cite{YaganMakowskiAllerton2009}.
In line with developments currently available for
other classes of graphs, e.g.,
Erd\H{o}s-R\'enyi graphs \cite[Chap. 3]{JansonLuczakRucinski}
and geometric random graphs \cite[Chap. 3]{PenroseBook},
it would be interesting to consider 
the containment problem for small subgraphs other than triangles
other than triangle in the context of random key graphs.
Given the difficulties encountered in the case of
This is likely to be a challenging problem
given the difficulties encountered in the simple case of
triangles.

The paper is organized as follows: In Section
\ref{sec:RandomKeyGraph} we formally introduce the class of random
key graphs while in Section \ref{sec:MainResults} we present the
main results of the paper given as Theorem \ref{thm:ZeroLaw} and
Theorem \ref{thm:OneLaw}.
Section \ref{sec:ComparingWithERGs} compares
these results with the corresponding zero-one law
in Erd\H{o}s-R\'enyi graphs.
The zero-one laws are established by an application
of the method of first and second moments, respectively 
\cite[p.  55]{JansonLuczakRucinski}. To that end, in Section
\ref{sec:FirstMoment}, we compute the expected value of the number
of triangles in random key graphs. Asymptotic results to be used
in the proofs of several results are then collected in Section
\ref{sec:UsefulAsymptotics} for easy reference. In Section
\ref{subsec:ProofTheoremZero}, we give a proof of the zero-law
(Theorem \ref{thm:ZeroLaw}) while an outline for the proof of the
one-law (Theorem \ref{thm:OneLaw}) is provided in Section
\ref{subsec:ProofTheoremOne}. The final sections of the paper,
namely Sections \ref{sec:SecondMoment} through
\ref{sec:FinalStep}, are devoted to completing the various steps
of the proof of Theorem \ref{thm:OneLaw}. Additional technical
derivations are given in Appendices \ref{Appendix:A},
\ref{Appendix:B} and \ref{Appendix:C}.

A word on the notation and conventions in use: All limiting
statements, including asymptotic equivalences, are understood with
$n$ going to infinity. The random variables (rvs) under
consideration are all defined on the same probability triple
$(\Omega, {\cal F}, \mathbb{P})$. Probabilistic statements are
made with respect to this probability measure $\mathbb{P}$, and we
denote the corresponding expectation operator by $\mathbb{E}$. The
indicator function of an event $E$ is denoted by $\1{E}$. For any
discrete set $S$ we write $|S|$ for its cardinality.

\section{Random key graphs}
\label{sec:RandomKeyGraph}

The model is parametrized by the number $n$ of nodes, the size $P$
of the key pool and the size $K$ of each key ring with $K \leq P$.
We often group the integers $P$ and $K$ into the ordered pair
$\theta \equiv (K,P)$ in order to simplify the notation. Now, for
each node $i=1, \ldots , n$, let $K_i(\theta)$ denote the random
set of $K$ distinct keys assigned to node $i$ and let
$\mathcal{P}$ be the set of all keys. 
The rvs $K_1(\theta), \ldots , K_n(\theta)$ are assumed 
to be {\em i.i.d.} rvs, each of which
is {\em uniformly} distributed with
\begin{equation}
\bP{ K_i(\theta) = S } = {P \choose K} ^{-1}, \qquad i=1,\ldots, n
\label{eq:KeyDistrbution1}
\end{equation}
for any subset $S$ of $\mathcal{P}$ which contains exactly $K$
elements. This corresponds to selecting keys randomly and {\em
without} replacement from the key pool.

Distinct nodes $i,j=1, \ldots , n$ are said to be adjacent if they
share at least one key in their key rings, namely
\begin{equation}
K_i (\theta) \cap K_j (\theta) \not = \emptyset ,
\label{eq:KeyGraphConnectivity}
\end{equation}
in which case an undirected link is assigned between nodes $i$ and
$j$. The resulting random graph defines the {\em random key graph}
on the vertex set $\{ 1, \ldots , n\}$, hereafter denoted
$\mathbb{K}(n; \theta )$.

For distinct $i,j =1, \ldots , n$, it is easy to check that
\begin{equation}
\bP{ K_i (\theta) \cap K_j (\theta) = \emptyset } = q(\theta)
\end{equation}
with
\begin{equation}
q (\theta) := \left \{
\begin{array}{ll}
0 & \mbox{if~ $P <2K$} \\
  &                 \\
\frac{{P-K \choose K}}{{P \choose K}} & \mbox{if~ $2K \leq P$,}
\end{array}
\right . \label{eq:q_theta}
\end{equation}
whence the probability of edge occurrence between any two nodes is
equal to $1-q(\theta)$. The expression given in (\ref{eq:q_theta})
is a simple consequence of the often used fact that
\begin{equation}
\bP{ S \cap K_i(\theta) = \emptyset } 
= 
\frac{{P- |S| \choose K}}{{P \choose K}}, 
\quad i=1, \ldots ,n
\label{eq:Probab_key_ring_does_not_intersect_S}
\end{equation}
for every subset $S$ of $\{ 1, \ldots , P \}$ with $|S| \leq P-K$.
Note that if $P<2K$ there exists an edge between any pair of
nodes, so that $\mathbb{K}(n;\theta)$ coincides with the complete
graph $K_{n}$. Also, we always have $0 \leq q(\theta) < 1 $ with
$q(\theta)> 0$ if and only if $2K \leq P$.

\section{The main results}
\label{sec:MainResults}

Pick positive integers $K$ and $P$ such that $K \leq P$. Fix $n=3,
4, \ldots $ and for distinct $i, j, k= 1, \ldots, n$, define the
indicator function
\[
\chi_{n,{ijk}}(\theta) := \1{ {\rm Nodes~}i, j ~{\rm and~} k~{\rm
form~a~triangle~in~}
                         \mathbb{K}(n; \theta) }.
\]
The number of (unlabelled) triangles in $\mathbb{K}(n; \theta)$ is
simply given by
\begin{equation}
T_n (\theta) := \sum_{(ijk)} \chi_{n,{ijk}}(\theta)
\label{eq:NumberOfTriangles}
\end{equation}
where $\sum_{(ijk)}$ denotes summation over all distinct triples
$ijk$ with $1 \leq i < j<k \leq n$. The event $T(n, \theta)$ that
there exists at least one triangle in $\mathbb{K}(n; \theta )$ is
then characterized by
\begin{equation}
T(n, \theta) := [ T_n (\theta) >  0 ] = [ T_n (\theta) = 0 ]^{c}.
\label{eq:triangle_basis}
\end{equation}

The main result of the paper is a zero-one law for the existence
of triangles in random key graphs. To state these results we find it
convenient to make use of the quantity
\begin{equation}
\tau(\theta)
:= \frac{K^3}{P^2} + \left ( \frac{K^2}{P} \right )^3,
\quad \theta = (K,P)
\label{eq:Tau(Theta)}
\end{equation}
with positive integers $K$ and $P$ such that $K \leq P$.
For simplicity of exposition we refer to any pair of functions
$P,K: \mathbb{N}_0 \rightarrow \mathbb{N}_0$ as a {\em scaling}
provided the natural conditions
\begin{equation}
K_n \leq P_n, \quad n=2, 3, \ldots 
\label{eq:ScalingCondition}
\end{equation}
are satisfied. The zero-law is given first.

\begin{theorem}
{\sl For any scaling $P,K: \mathbb{N}_0 \rightarrow \mathbb{N}_0$,
we have the zero-law
\begin{equation}
\lim_{n \rightarrow \infty } \bP{T(n, \theta_n)} = 0
\label{eq:MainTheoremZero}
\end{equation}
under the condition
\begin{equation}
\lim_{n \rightarrow \infty } n^3 \tau(\theta_n) = 0 .
\label{eq:ConditionForZero}
\end{equation}
} \label{thm:ZeroLaw}
\end{theorem}

The one-law given next assumes a more involved form.

\begin{theorem}
{\sl For any scaling $P,K: \mathbb{N}_0 \rightarrow \mathbb{N}_0$
for which the limit $\lim_{n \rightarrow \infty } q(\theta_n) = q ^\star $ 
exists, we have the one-law
\begin{equation}
\lim_{n \rightarrow \infty } \bP{T(n, \theta_n)} = 1
\label{eq:MainTheoremOne}
\end{equation}
either if $ 0 \leq q^\star < 1$ or if $q^\star =1$ under the
condition
\begin{equation}
\lim_{n \rightarrow \infty } n^3 \tau(\theta_n) = \infty .
\label{eq:ConditionForOne}
\end{equation}
} \label{thm:OneLaw}
\end{theorem}

Theorem \ref{thm:ZeroLaw} and Theorem \ref{thm:OneLaw} will be
established by the method of first and second moments,
respectively \cite[p. 55]{JansonLuczakRucinski}, applied to
the count variables defined at (\ref{eq:NumberOfTriangles}). 
To facilitate comparison with Erd\H{o}s-R\'enyi graphs, we combine
Theorem \ref{thm:ZeroLaw} and Theorem \ref{thm:OneLaw}
into a symmetric, but somewhat weaker, statement.

\begin{theorem}
{\sl For any scaling $P,K: \mathbb{N}_0 \rightarrow \mathbb{N}_0$
for which $\lim_{n \rightarrow \infty } q(\theta_n) = 1$, we have
\begin{eqnarray}
\lim_{n \rightarrow \infty } \bP{ T(n; \theta_n ) }
 = \left \{
\begin{array}{ll}
0 & \mbox{if~ $\lim_{ n\rightarrow \infty } n^3  \tau (\theta_n) = 0 $} \\
  &                 \\
1 & \mbox{if~ $\lim_{ n\rightarrow \infty } n^3  \tau (\theta_n) =
\infty $.}
\end{array}
\right . 
\label{eq:TriangleZeroOneLaw+RKG}
\end{eqnarray}
} 
\label{thm:TriangleZeroOneLaw+RKG}
\end{theorem}

\section{Comparing with Erd\H{o}s-R\'enyi graphs}
\label{sec:ComparingWithERGs}

In this section we compare
Theorem \ref{thm:TriangleZeroOneLaw+RKG} 
with its analog for Erd\H{o}s-R\'enyi graphs.
First some notation:
For each $p$ in $[0,1]$ and $n=2,3, \ldots $, let
$\mathbb{G}(n;p)$ denote the Erd\H{o}s-R\'enyi graph on the vertex
set $\{ 1, \ldots , n\}$ with edge probability $p$. 
In analogy
with (\ref{eq:NumberOfTriangles}) and (\ref{eq:triangle_basis})
let $T_n(p)$ denote
the number of (unlabelled) triangles in $\mathbb{G}(n; p)$, and
define $T(n,p)$ as the event that
there exists at least one triangle in $\mathbb{G}(n; p)$, i.e.,
$T(n,p) = [ T_n(p) > 0 ]$.
we also refer to any mapping
$p: \mathbb{N}_0 \rightarrow [0,1]$ as a scaling for
Erd\H{o}s-R\'enyi graphs. The following zero-one law for
connectivity in Erd\H{o}s-R\'enyi graphs is well known \cite{ER1960}.

\begin{theorem}
{\sl For any scaling $p: \mathbb{N}_0 \rightarrow [0,1]$, we have
\begin{eqnarray}
\lim_{n \rightarrow \infty } \bP{ T(n; p_n ) }
 = \left \{
\begin{array}{ll}
0 & \mbox{if~ $\lim_{ n\rightarrow \infty } n^3  \tau ^ \star (p_n) = 0 $} \\
  &                 \\
1 & \mbox{if~ $\lim_{ n\rightarrow \infty } n^3  \tau ^ \star (p_n) = \infty $}
\end{array}
\right . \label{eq:ERZeroOneLawTriangle}
\end{eqnarray}
where
\begin{equation}
\tau ^ \star (p) := p^3, \quad p \in [0,1] .
\end{equation}
}
\label{thm:ERZeroOneLawTriangle}
\end{theorem}

As this result is also established by the method of first and second
moments, its form is easily understood once we note that
\begin{equation}
\bE{ T_n(p) } = {n \choose 3} \tau ^ \star (p),
\quad 0 \leq p \leq 1
\label{eq:ER+FirstMoment}
\end{equation}
for all $n=3,4,  \ldots$.

As mentioned earlier, random key graphs are {\em not} equivalent
to Erd\H{o}s-Renyi graphs even when their edge probabilities are
matched, i.e., $\mathbb{G}(n;p) \neq_{st} \mathbb{K}(n; \theta)$
with $p = 1 - q(\theta)$; see \cite{YaganMakowskiAllerton2009}
for a discussion of similarities.
However, in order to meaningfully compare the zero-one law of Theorem
\ref{thm:ERZeroOneLawTriangle} with that contained in Theorem
\ref{thm:TriangleZeroOneLaw+RKG},
we say that the scaling $p: \mathbb{N}_0 \rightarrow [0,1]$
(for Erd\H{o}s-R\'enyi graphs) is {\em asymptotically matched}
to the scaling $P,K: \mathbb{N}_0 \rightarrow \mathbb{N}_0$
(for random key graphs) if
\begin{equation}
p_n \sim 1 - q(\theta_n).
\label{eq:AsymptoticMatchingCondition}
\end{equation}
This is equivalent to requiring that the expected average degrees
are asymptotically equivalent. Under the
natural condition $\lim_{n \rightarrow \infty} q(\theta_n) = 1$,
the matching condition (\ref{eq:AsymptoticMatchingCondition})
amounts to
\begin{equation}
p_n \sim \frac{K_n^2}{P_n} \label{eq:AsymptoticMatchingEquiv}
\end{equation}
by virtue of Lemma \ref{lem:AsymptoticEquivalence1}.

The definitions readily yield
\[
\frac{\tau(\theta_n)}
     {\tau^\star(p_n)}
= \frac{1}{p_n^3} \cdot \left ( \frac{K^3_n}{P^2_n} \right )
+
\frac{1}{p_n^3} \cdot \left ( \frac{K^2_n}{P_n} \right )^3,
\quad n=2,3, \ldots
\]
whence
\begin{equation}
\frac{\tau(\theta_n)}{\tau ^ \star (p_n)}
\sim
1 + \frac{P_n}{K_n^3}
\label{eq:AsymptoticRatio}
\end{equation}
under (\ref{eq:AsymptoticMatchingEquiv}).
By Proposition \ref{prop:AsymptoticEquivalence2},
this last statement is equivalent to
\begin{equation}
\frac{ \bE{ T_n(\theta_n) } }
     { \bE{ T_n(p_n) } }
\sim
1 + \frac{P_n}{K_n^3}
\label{eq:FirstMomentAsymptoticRatio}
\end{equation}
as we make use of the expressions (\ref{eq:ER+FirstMoment})
and (\ref{eq:FirstMoment}).
In other words, for large $n$ the expected number of triangles
in random key graphs is always at least as large as
the corresponding quantity
in asymptotically matched Erd\H{o}s-R\'enyi graphs.

In the context of WSNs, it is natural to select the parameters
$K_n$ and $P_n$ of the scheme of Eschenauer and Gligor such that
the induced random key graph is {\em connected}. However, there is
a tradeoff between connectivity and security
\cite{DiPietroManciniMeiPanconesiRadhakrishnan2008}. This requires
that $\frac{K_n ^ 2}{P_n}$ be kept as close as possible to the
critical scaling $\frac{\log n}{n}$ for connectivity;
see the papers \cite{BlackburnGerke,
DiPietroManciniMeiPanconesiRadhakrishnan2008, Rybarczyk2009,
YaganMakowskiISIT2009, YaganMakowskiConnectivity}.
In the desired near boundary regime, this amounts to
\begin{equation}
\frac{K_n ^ 2 }{ P_n } \sim c \cdot \frac{ \log n }{ n }
\label{eq:K^2/P_sim_logn/n}
\end{equation}
with $c > 1$ but close to one, and from
(\ref{eq:FirstMomentAsymptoticRatio}) we see that
\begin{equation}
\frac{ \bE{ T_n(\theta_n) } }
     { \bE{ T_n(p_n) } }
\sim 1 
\quad \mbox{if and only if} 
\quad K_n \gg \frac{n}{\log n} . 
\label{eq:NotPractical}
\end{equation}
The expected number of triangles in random key graphs is then of
the same order as the corresponding quantity in asymptotically
matched Erd\H{o}s-R\'enyi graphs with $ \bE{ T_n(\theta_n) } \sim
\bE{ T_n(p_n) } \sim \frac{c^3}{6} \left ( \log n \right )^3$.
This conclusion holds regardless of the value of $c$ in
(\ref{eq:K^2/P_sim_logn/n}).

However, given the limited memory and computational power of the sensor
nodes, the key ring sizes at (\ref{eq:NotPractical}) are not
practical. In addition, they will lead to {\em high} node degrees
and this in turn will decrease network {\em resiliency} against
node capture attacks. Indeed, in
\cite[Thm. 5.3]{DiPietroManciniMeiPanconesiRadhakrishnan2008} it was
proposed that security in WSNs be ensured by selecting $K_n$ and $P_n$
such that $\frac{K_n}{P_n} \sim \frac{1}{n}$, a requirement which
then leads to
\begin{equation}
K_n \sim c \cdot \log n \label{eq:K_n_sim_logn}
\end{equation}
under (\ref{eq:K^2/P_sim_logn/n}), and
(\ref{eq:FirstMomentAsymptoticRatio}) implies
\begin{equation}
\lim_{n \rightarrow \infty}
\frac{ \bE{ T_n(\theta_n) } }
     { \bE{ T_n(p_n) } }
=
\lim_{n \rightarrow \infty}
\left ( 1 + \frac{n}{(c \cdot \log n) ^2} \right )
= \infty.
\end{equation}
Hence, for realistic WSN scenarios the expected number of triangles
in the induced random key graphs can be orders of magnitude larger
than in Erd\H{o}s-R\'enyi graphs.
This provides a clear example where transferring known
results for Erd\H{o}s-R\'enyi graphs to random key graphs by
asymptotically matching their edge probabilities can be misleading.

\section{Computing the first moment}
\label{sec:FirstMoment}

With positive integers $K$ and $P$ such that $K \leq P$, define
\begin{equation}
\beta(\theta) := (1-q(\theta))^3 + q(\theta)^3 - q(\theta)
r(\theta) \label{eq:beta(tetha)}
\end{equation}
where we have set
\begin{equation}
r (\theta) := \left \{
\begin{array}{ll}
0 & \mbox{if~ $P <3K$} \\
  &                 \\
\frac{{P-2K \choose K}}{{P \choose K}} & \mbox{if~ $3K \leq P$.}
\end{array}
\right . \label{eq:r(tetha)}
\end{equation}
Direct inspection shows that
\begin{equation}
r(\theta) \leq q(\theta)^2 \label{eq:r(tetha)B}
\end{equation}
whence
\begin{equation}
\beta (\theta) \geq (1 - q(\theta) )^3 > 0 . \label{eq:r(tetha)C}
\end{equation}

\begin{lemma}
{\sl For positive integers $K$ and $P$ such that $K \leq P$, we
have
\begin{equation}
\bE{ T_n (\theta) } = {n \choose 3 } \beta ( \theta ), \quad
n=3,4, \ldots \label{eq:FirstMoment}
\end{equation}
} \label{lem:FirstMoment}
\end{lemma}

To help deriving (\ref{eq:FirstMoment}) we introduce the events
\begin{equation}
A(\theta) := [ K_1 (\theta) \cap K_2 (\theta) \neq \emptyset ]
\cap [ K_1 (\theta) \cap K_3 (\theta) \neq \emptyset ]
\end{equation}
and
\begin{eqnarray}
B(\theta) &:= & [ K_1 (\theta) \cap K_2 (\theta) \neq \emptyset ]
\cap [ K_1 (\theta) \cap K_3 (\theta) \neq \emptyset ] \cap [ K_2
(\theta) \cap K_3 (\theta) \neq \emptyset ]
\nonumber \\
&=& A(\theta) \cap [ K_2 (\theta) \cap K_3 (\theta) \neq \emptyset
].
\end{eqnarray}
The event $A(\theta)$ captures the existence of edges between node
$1$ and the pair of nodes $2$ and $3$, respectively, in
$\mathbb{K}(n;\theta)$, while $B(\theta)$ is the event where the
nodes $1$, $2$ and $3$ form a triangle in $\mathbb{K}(n;\theta)$.

\begin{lemma}
The probability of the event $A(\theta)$ is given by
\begin{equation}
\bP{ A(\theta) } = (1-q(\theta))^2. \label{eq:A_theta}
\label{eq:A_theta}
\end{equation}
\label{lem:A_theta}
\end{lemma}

In the proof of Lemma \ref{lem:A_theta} (as well as in other
proofs) we omit the explicit dependence on $\theta$ when no
confusion arises from doing so.

\myproof Under the enforced independence assumptions we note that
\begin{eqnarray}
\bP{ A(\theta) } &=& \sum_{ |S| = K } \bP{ K_1 = S,
                      S \cap K_2 \neq \emptyset,
                      S \cap K_3 \neq \emptyset }
\nonumber \\
&=& \sum_{ |S| = K } \bP{ K_1 = S}
                      \bP{ S \cap K_2 \neq \emptyset }
                      \bP{ S \cap K_3 \neq \emptyset }
\nonumber \\
&=& \left ( 1 - q(\theta) \right )^2
\end{eqnarray}
as we make use of (\ref{eq:Probab_key_ring_does_not_intersect_S})
with $\sum_{ |S| = K } \bP{ K_1 = S} = 1 $. \myendpf

In many of the forthcoming calculations we make repeated use of
the fact that for any pair of events, say $E$ and $F$, we have
\begin{equation}
\bP{ E \cap F } = \bP{ E} - \bP{ E \cap F^c } .
\label{eq:Set_Formula}
\end{equation}
In particular, we can now conclude from Lemma \ref{lem:A_theta}
that
\begin{eqnarray}
\lefteqn{\bP{ K_1(\theta)\cap K_2(\theta) = \emptyset,\;
     K_1(\theta) \cap K_3(\theta) \neq \emptyset }} &&
\nonumber \\
&=& \bP{ K_1(\theta)\cap K_2(\theta) \neq \emptyset,\;
     K_1(\theta) \cap K_3(\theta) = \emptyset }
\nonumber \\
&=& q(\theta)(1-q(\theta)) \label{eq:cor_q_1-q}
\end{eqnarray}
and
\begin{equation}
\bP{ K_1(\theta)\cap K_2(\theta) = \emptyset,\;
     K_1(\theta) \cap K_3(\theta) = \emptyset }
= q(\theta)^2 . \label{eq:cor_q^2}
\end{equation}
These facts will now be used in computing the probability of
$B(\theta)$.

\begin{lemma}
{\sl With $\beta (\theta)$ given at (\ref{eq:beta(tetha)}) we have
\begin{equation}
\bP{ B(\theta) } = \beta (\theta) . \label{eq:B_theta}
\end{equation}
} \label{lem:Prob_B}
\end{lemma}

\myproof Repeated use of (\ref{eq:Set_Formula}) yields
\begin{eqnarray}
\bP{ B(\theta) } &=& \bP{ K_1\cap K_2\neq\emptyset,\; K_1 \cap K_3
\neq\emptyset }
\nonumber \\
& & ~ - \bP{ K_1\cap K_2\neq\emptyset,\;
             K_1 \cap K_3 \neq\emptyset,\;
             K_2 \cap K_3 =\emptyset }
\nonumber \\
&=& \bP{ A(\theta) } - \bP{ K_1 \cap K_2 \neq\emptyset,\;
                            K_2 \cap K_3 =\emptyset }
\nonumber \\
& & ~ + \bP{ K_1 \cap K_2 \neq \emptyset,\;
             K_1 \cap K_3 = \emptyset,\;
             K_2 \cap K_3 =\emptyset }
\nonumber \\
&=& (1-q(\theta))^2 - q(\theta)(1-q(\theta))
    + \bP{ K_1 \cap K_3 = \emptyset,\;
           K_2 \cap K_3 =\emptyset }
\nonumber \\
& & ~ - \bP{ K_1\cap K_2=\emptyset,\; K_1 \cap K_3 = \emptyset,\;
             K_2 \cap K_3 =\emptyset }
\nonumber \\
&=& (1-q(\theta))^2 - q(\theta)(1-q(\theta)) + q(\theta)^2
\nonumber \\
& & ~ - \bP{ K_1\cap K_2=\emptyset,\; K_1 \cap K_3 = \emptyset,\;
             K_2 \cap K_3 =\emptyset }
\end{eqnarray}
as we recall (\ref{eq:A_theta}), (\ref{eq:cor_q_1-q}) and
(\ref{eq:cor_q^2}).

By independence we get
\begin{eqnarray}
\lefteqn{ \bP{ K_1\cap K_2=\emptyset,\; K_1 \cap K_3 =
\emptyset,\;
             K_2 \cap K_3 =\emptyset }
} & &
\nonumber \\
&=& \bP{ K_1\cap K_2=\emptyset,\;
             (K_1 \cup K_2 ) \cap K_3 =\emptyset }
\nonumber \\
&=& \sum_{|S|=|T| = K, S \cap T = \emptyset} \bP{ K_1 = S, K_2 = T
} \bP{ (S \cup T ) \cap K_3 =\emptyset }
\nonumber \\
&=& \sum_{|S|=|T|= K, S \cap T = \emptyset} \bP{ K_1 = S , K_2 = T
} \cdot r(\theta )
\nonumber \\
&=& \bP{ K_1 \cap K_2 = \emptyset } \cdot r(\theta )
\end{eqnarray}
by invoking (\ref{eq:Probab_key_ring_does_not_intersect_S}) (since
$|S \cup T| = 2K$ under the constraints $|S|=|T|= K$ and $S \cap T
= \emptyset$). Thus,
\[
\bP{ B(\theta) } = (1-q(\theta))^2 - q(\theta)(1-q(\theta)) +
q(\theta)^2 - q(\theta) r(\theta),
\]
and the desired result follows upon noting the relation
\[
(1-q(\theta))^2 - q(\theta)(1-q(\theta)) + q(\theta)^2 =
(1-q(\theta))^3 + q(\theta)^3 .
\]
\myendpf

The proof of Lemma \ref{lem:FirstMoment} is now straightforward:
Fix $n=3,4, \ldots $. Exchangeability yields
\begin{equation}
\bE{ T_n (\theta)} = {n \choose 3} \bE{ \chi_{n,{123}} (\theta) }
\label{eq:FirstMomentExpression}
\end{equation}
and the desired conclusion follows as we make use of Lemma
\ref{lem:Prob_B}.

\section{Some useful asymptotics}
\label{sec:UsefulAsymptotics}

In this section we collect a number of asymptotic results that
prove useful in establishing some of the results derived in this
paper. The first result was already obtained in
\cite{YaganMakowskiConnectivity}.

\begin{lemma}
{\sl For any scaling $P,K: \mathbb{N}_0 \rightarrow \mathbb{N}_0$,
we have
\begin{equation}
\lim_{n \rightarrow \infty} q(\theta_n) = 1 \label{eq:Condition1}
\end{equation}
if and only if
\begin{equation}
\lim_{n \rightarrow \infty} \frac{K^2_n}{P_n} = 0,
\label{eq:Condition2}
\end{equation}
and under either condition the asymptotic equivalence
\begin{equation}
1 - q(\theta_n) \sim \frac{K^2_n}{P_n}
\label{eq:AsymptoticsEquivalence1}
\end{equation}
holds.
} 
\label{lem:AsymptoticEquivalence1}
\end{lemma}

Since $1 \leq K_n \leq {K_n}^2$ for all $n=1,2, \ldots $, the
condition (\ref{eq:Condition2}) implies
\begin{equation}
\lim_{n \rightarrow \infty } \frac{K_n}{P_n} = 0
\label{eq:RatioConditionStrong+Consequence1}
\end{equation}
and
\begin{equation}
\lim_{n \rightarrow \infty } P_n = \infty .
\label{eq:RatioConditionStrong+Consequence2}
\end{equation}
so that for any $c > 0$, we have
\begin{equation}
c K_n < P_n \label{eq:Condition0}
\end{equation}
for all $n$ sufficiently large in $\mathbb{N}_0$ (dependent on
$c$).

The following asymptotic equivalence will be crucial to stating
the results in a more explicit form.

\begin{proposition}
{\sl For any scaling $P,K: \mathbb{N}_0 \rightarrow \mathbb{N}_0$
satisfying (\ref{eq:Condition1})-(\ref{eq:Condition2}), we have
the asymptotic equivalence
\begin{equation}
\beta(\theta_n) \sim \tau(\theta_n) .
\label{eq:AsymptoticsEquivalence2}
\end{equation}
}
\label{prop:AsymptoticEquivalence2}
\end{proposition}

\myproof From (\ref{eq:beta(tetha)}), we get
\[
\beta ( \theta_n ) = \left ( 1 - q(\theta_n) \right )^3
                        +  q ( \theta_n )^3
\left( 1-\frac { r ( \theta_n ) } { q ^ 2 (\theta_n) } \right) .
\]
Under the enforced assumptions Lemma
\ref{lem:AsymptoticEquivalence1} already implies
\[
\left ( 1 - q(\theta_n) \right )^3 \sim \left ( \frac{K^2_n}{P_n}
\right )^3
\]
with $q ( \theta_n )^ 3 \sim 1$. It is now plain that the
equivalence (\ref{eq:AsymptoticsEquivalence2}) will hold if we
show that
\begin{equation}
1 - \frac{r(\theta_n)}{q(\theta_n)^2} \sim \frac{K^3_n}{P^2_n} .
\label{eq:AsymptoticsEquivalence2Reduced}
\end{equation}
This key technical fact is established in Appendix
\ref{Appendix:A}. \myendpf

The final result of this section also relies on Lemma
\ref{lem:AsymptoticEquivalence1}, and will prove useful in
establishing the one-law.

\begin{proposition}
{\sl  For any scaling $P,K: \mathbb{N}_0\rightarrow\mathbb{N}_0$
satisfying (\ref{eq:Condition1})-(\ref{eq:Condition2}), we have
\begin{equation}
\lim_{n \to \infty} n^2 (1-q(\theta_n))=\infty
\label{eq:n^2_1_q_to_inf}
\end{equation}
provided the condition (\ref{eq:ConditionForOne}) holds. }
\label{prop:n^2_1_q_to_inf}
\end{proposition}

\myproof Consider a scaling $P,K:
\mathbb{N}_0\rightarrow\mathbb{N}_0$ satisfying
(\ref{eq:Condition1})-(\ref{eq:Condition2}). By Lemma
\ref{lem:AsymptoticEquivalence1} the desired conclusion
(\ref{eq:n^2_1_q_to_inf}) will be established if we show
\begin{equation}
\lim_{n \to \infty} n^2 \frac{K^2_n}{P_n} = \infty .
\label{eq:n^2_1_q_to_infB}
\end{equation}
As condition (\ref{eq:ConditionForOne}) reads
\[
\lim_{n \to \infty} n^3 \left( \frac{K_n^3}{P_n^2} +
\left(\frac{K_n^2}{P_n}\right)^3 \right) = \infty ,
\]
we immediately get (\ref{eq:n^2_1_q_to_infB}) from it by virtue of
the trivial bounds
\[
n^3 \left(\frac{K_n^2}{P_n}\right)^3 = \left(\frac{n
K_n^2}{P_n}\right)^3 \leq \left(\frac{n^2 K_n^2}{P_n}\right)^3
\]
and
\[
n^3 \frac{K_n^3}{P_n^2} \leq n^4 \frac{K_n^4}{P_n^2} =
\left(\frac{n^2 K_n^2}{P_n}\right)^2
\]
valid for all $n=1,2, \ldots $. \myendpf

Proposition \ref{prop:n^2_1_q_to_inf} will be used as follows:
Pick $a>0$ and $b>0$, and consider a scaling $P,K:
\mathbb{N}_0\rightarrow\mathbb{N}_0$ satisfying
(\ref{eq:Condition1})-(\ref{eq:Condition2}). For each $n=2,3,
\ldots $, we get
\begin{eqnarray}
\frac{1}{n^2} \cdot \frac{ \left ( 1 - q (\theta_n) \right )^a }{
\beta (\theta_n)^b } &\leq& \frac{1}{n^2} \cdot \frac{ \left ( 1 -
q (\theta_n) \right )^a }
     { \left ( 1 - q (\theta_n) \right )^{3b} }
\nonumber \\
&=& \frac{1}{n^2 \left ( 1 - q (\theta_n) \right )} \cdot { \left
( 1 - q (\theta_n) \right )^{a-3b+1} } .
\end{eqnarray}
Therefore, under condition (\ref{eq:ConditionForOne}) Proposition
\ref{prop:n^2_1_q_to_inf} yields
\begin{equation}
\lim_{n \to \infty} \frac{1}{n^2} \cdot \frac{ \left ( 1 - q
(\theta_n) \right )^a }{ \beta (\theta_n)^b } = 0 \quad \mbox{if~}
a-3b+1 \geq 0 \label{eq:AsymptoticsForRatio}
\end{equation}
as we make use of (\ref{eq:Condition1})-(\ref{eq:Condition2}).

\section{Proofs of Theorem \ref{thm:ZeroLaw} and Theorem \ref{thm:OneLaw}}
\label{sec:ProofsTheoremsZeroOne}

\subsection{A proof of Theorem \ref{thm:ZeroLaw}}
\label{subsec:ProofTheoremZero}

Fix $n=3,4, \ldots $, An elementary bound for $\mathbb{N}$-valued
rvs yields
\begin{equation}
\bP{ T_n (\theta_n) > 0 } \leq \bE{ T_n(\theta_n) } ,
\end{equation}
so that
\begin{equation}
\bP{ T (n, \theta_n) } \leq {n \choose 3} \beta (\theta_n) .
\end{equation}
The conclusion (\ref{eq:MainTheoremZero}) follows if we show that
\begin{equation}
\lim_{n \rightarrow \infty} {n \choose 3} \beta (\theta_n) = 0
\label{eq:ConditionForZero+A}
\end{equation}
under (\ref{eq:ConditionForZero}).

The condition $\lim_{n \rightarrow \infty} n^3 \tau(\theta_n) = 0$
implies $\lim_{n \rightarrow \infty} \tau(\theta_n) = 0$ and
(\ref{eq:Condition2}) automatically holds. By Proposition
\ref{prop:AsymptoticEquivalence2} we conclude $\beta(\theta_n)
\sim \tau(\theta_n)$, whence $n^3 \beta(\theta_n) \sim n^3
\tau(\theta_n)$, and condition (\ref{eq:ConditionForZero}) is
indeed equivalent to (\ref{eq:ConditionForZero+A}) 
since ${n \choose 3} \sim \frac{n^3}{6}$.

\subsection{A proof of Theorem \ref{thm:OneLaw}}
\label{subsec:ProofTheoremOne}

Assume first that $q^\star$ satisfies $0 \leq q^\star < 1$. Fix
$n=3,4, \ldots $ and partition the $n$ nodes into the $k_n+1$
non-overlapping groups $(1,2,3)$, $(4,5,6)$, $\ldots $,
$(3k_n+1,3k_n+2,3k_n+3)$ with $k_n = \lfloor \frac{n-3}{3} \rfloor
$. If $\mathbb{K}(n;\theta_n)$ contains no triangle, then {\em
none} of these $k_n + 1$ groups of nodes forms a triangle. 
With this in mind we get
\begin{eqnarray}
\lefteqn{\bP{ T_n(\theta_n) = 0 } } & &
\nonumber \\
&\leq& \bP{ \bigcap_{\ell=0}^{k_n} \left [
\begin{array}{c}
\mbox{Nodes $3\ell+1,3\ell+2, 3\ell+3$ do not form } \\
\mbox{a triangle in $\mathbb{K}(n;\theta_n)$ } \\
\end{array}
\right ] }
\nonumber \\
&=& \prod_{\ell=0}^{k_n} \bP{
\begin{array}{c}
\mbox{Nodes $3\ell+1,3\ell+2, 3\ell+3$ do not form } \\
\mbox{a triangle in $\mathbb{K}(n;\theta_n)$ } \\
\end{array}
}
\label{eq:IndependenceTriangle} \\
&=& \left ( 1 - \beta(\theta_n) \right )^{k_n+1}
\nonumber \\
&\leq& \left ( 1 - (1-q(\theta_n) )^3 \right )^{k_n+1}
\label{eq:OneLawInequality0}
\\
&\leq& e^{- (k_n +1 ) (1-q(\theta_n) )^3 }.
\label{eq:OneLawInequality1}
\end{eqnarray}
Note that (\ref{eq:IndependenceTriangle}) follows from the fact
that the events
\[
\left [
\begin{array}{c}
\mbox{Nodes $3\ell+1,3\ell+2, 3\ell+3$ do not form } \\
\mbox{a triangle in $\mathbb{K}(n;\theta_n)$ } \\
\end{array}
\right ], \quad \ell =0, \ldots , k_n
\]
are mutually independent due to the non-overlap condition, while
the inequality (\ref{eq:OneLawInequality0}) is justified with the
help of (\ref{eq:r(tetha)C}). Let $n$ go to infinity in the
inequality (\ref{eq:OneLawInequality1}). 
The condition $q^\star < 1$ implies
$\lim_{n \rightarrow \infty} \bP{ T(n,\theta_n )^c } = 0$ 
since $k_n \sim \frac{n}{3}$ so that
$\lim_{n \rightarrow \infty} ( k_n + 1 ) (1-q(\theta_n) )^3 = \infty $. 
This establishes (\ref{eq:MainTheoremOne}).

To handle the case $q^\star =1$, we use a standard bound which
forms the basis of the method of second moment 
\cite[remark 3.1, p. 55]{JansonLuczakRucinski}. Here it takes the form
\begin{equation}
\frac{ \bE{ T_n (\theta_n)}^2}{ \bE{T_n (\theta_n)^2 }} \leq
\bP{T_n (\theta_n) > 0 }, \quad n=3,4, \ldots
\label{eq:SecondMoment+b}
\end{equation}
It is now plain that (\ref{eq:MainTheoremOne}) will be established
in the case $q^\star =1 $ if we show the following result.

\begin{proposition}
{\sl For any scaling $P,K: \mathbb{N}_0 \rightarrow \mathbb{N}_0$
satisfying (\ref{eq:Condition1})-(\ref{eq:Condition2}), we have
\begin{equation}
\lim_{n \rightarrow \infty } 
\frac{\bE{T_n (\theta_n)^2}} {\bE{T_n (\theta_n)}^2} 
=1 
\label{eq:OneLawConvergenceSecondMoment}
\end{equation}
under the condition (\ref{eq:ConditionForOne}). }
\label{prop:OneLawConvergenceSecondMoment}
\end{proposition}

The remainder of the paper is devoted to establishing Proposition
\ref{prop:OneLawConvergenceSecondMoment}. As will soon become
apparent this is a bit quite more involved than expected.

\section{Computing the second moment}
\label{sec:SecondMoment}

A natural step towards establishing Proposition
\ref{prop:OneLawConvergenceSecondMoment} consists in computing the
second moment of the count variables (\ref{eq:NumberOfTriangles}).

\begin{proposition}
{\sl For positive integers $K$ and $P$ such that $K \leq P$, we
have
\begin{eqnarray}
\bE{ T_n (\theta)^2 } = \bE{ T_n (\theta) } 
&+& \left ( \frac{ {n-3\choose 3} }{ {n \choose 3} } 
+ 3 \frac{ {n-3\choose 2} }{ {n \choose 3} } \right ) 
\cdot \bE{ T_n (\theta) }^2
\label{eq:SecondMoment} \\
&+& {n \choose 3}{3 \choose 2}{n-3\choose 1} \cdot
\bE{\chi_{n,123}(\theta) \chi_{n, 124} (\theta)} \nonumber
\end{eqnarray}
for all $n=3,4, \ldots$ with
\begin{eqnarray}
\lefteqn{ \bE{\chi_{n,123}(\theta) \chi_{n, 124} (\theta)} } & &
\nonumber \\
&=& - (1-q(\theta))^5 + 2 \left ( 1- q(\theta) \right )^2 \beta
(\theta)
\nonumber \\
& & ~ - \frac{1}{q(\theta)} \left ( \beta (\theta) -
(1-q(\theta))^3 \right )^2
      + \sum_{k=0}^K c_k (\theta) - q(\theta)^4
\label{eq:SecondMomentCross}
\end{eqnarray}
where we have set
\begin{equation}
c_k(\theta) := {{\begingroup {K} \choose {k} \endgroup}
 {\begingroup {P-K} \choose {K-k} \endgroup}
\over{P \choose K}} \cdot \left(  {{\begingroup {P-2K+k} \choose
{K} \endgroup}
       \over{P \choose K}}
\right)^2, \quad k=0,1, \ldots , K. \label{eq:SecondMomentCrossB}
\end{equation}
} \label{prop:SecondMoment}
\end{proposition}

A careful inspection of the definition (\ref{eq:subs5to6}) 
given for the quantities (\ref{eq:SecondMomentCrossB}) yields the
probabilistic interpretation
\begin{equation}
c_k (\theta) = \bP{ |K_1(\theta) \cap K_2 (\theta) | = k,
        \left ( K_1(\theta) \cup K_2(\theta) \right ) \cap K_i(\theta)
       = \emptyset, \ i=3,4 }
\label{eq:SecondMomentCrossBInterpretation}
\end{equation}
for each $k=0,1, \ldots , K$.

\myproof Consider positive integers $K$ and $P$ such that $K \leq
P$ and fix $n=3, 4, \ldots $. By exchangeability and by the binary
nature of the rvs involved we readily conclude that
\begin{eqnarray}
\bE{ T_n (\theta)^2 } &=& \sum_{(ijk)}\sum_{(abc)}
\bE{\chi_{n,{ijk}} (\theta) \chi_{n,{abc}} (\theta) }
\nonumber \\
&=& \bE{ T_n (\theta)} \label{eq:SecondMomentExpression}
\nonumber \\
& & +  {n \choose 3}{3 \choose 2}{n-3\choose 1}
        \bE{\chi_{n,123}(\theta) \chi_{n, 124} (\theta)}
\nonumber \\
& & +  {n \choose 3}{3 \choose 1}{n-3\choose 2}
       \bE{\chi_{n,123}(\theta) \chi_{n, 145} (\theta)}
\nonumber \\
& & +  {n \choose 3}{n-3\choose 3}
        \bE{\chi_{n,123}(\theta) \chi_{n, 456} (\theta)}.
\label{eq:SecondMomentPart0}
\end{eqnarray}

Under the enforced independence assumptions the rvs
$\chi_{n,123}(\theta)$ and $\chi_{n, 456} (\theta)$ are
independent and identically distributed. As a result,
\[
\bE{\chi_{n,123}(\theta) \chi_{n, 456} (\theta)} =
\bE{\chi_{n,123}(\theta)} \bE{ \chi_{n, 456} (\theta)} = \beta
(\theta)^2
\]
so that
\begin{equation}
{n \choose 3}{n-3\choose 3} \bE{\chi_{n,123}(\theta) \chi_{n, 456}
(\theta)} = \frac{ {n-3\choose 3} }{ {n \choose 3} } \cdot \bE{
T_n (\theta) }^2 \label{eq:SecondMomentPart1}
\end{equation}
as we make use of the relation (\ref{eq:FirstMoment}).

On the other hand, we readily check that the indicator rvs
$\chi_{n,123}(\theta)$ and $\chi_{n,145}(\theta)$ are independent
and identically distributed {\em conditionally} on $K_1(\theta)$
with
\[
\bP{ \chi_{n,123}(\theta) = 1 | K_1 (\theta ) = S } = \bP{
\chi_{n,123}(\theta) = 1 } = \beta (\theta), \quad S \in {\cal
P}_K.
\]
A similar statement applies to $\chi_{n,145}(\theta)$, and the rvs
$\chi_{n,123}(\theta)$ and $\chi_{n,145}(\theta)$ are therefore
(unconditionally) independent and identically distributed so that
\[
\bE{\chi_{n,123}(\theta) \chi_{n, 145} (\theta)} =
\bE{\chi_{n,123}(\theta) } \bE{ \chi_{n, 145} (\theta)}.
\]
As before this last observation yields
\begin{equation}
{n \choose 3}{3 \choose 1}{n-3\choose 2} \bE{\chi_{n,123}(\theta)
\chi_{n, 145} (\theta)} = 3 \frac{ {n-3\choose 2} }{ {n \choose 3}
} \cdot \bE{ T_n (\theta) }^2 \label{eq:SecondMomentPart2}
\end{equation}
by virtue of (\ref{eq:FirstMoment}).

The evaluation
(\ref{eq:SecondMomentCross})--(\ref{eq:SecondMomentCrossB}) of the
moment $\bE{\chi_{n,123}(\theta) \chi_{n, 124} (\theta)}$ is
rather lengthy, although quite straightforward; details are given
in Appendix \ref{Appendix:B}. Reporting
(\ref{eq:SecondMomentCross})--(\ref{eq:SecondMomentCrossB}),
(\ref{eq:SecondMomentPart1}) and (\ref{eq:SecondMomentPart2}) into
(\ref{eq:SecondMomentPart0}) establishes Proposition
\ref{prop:SecondMoment}. \myendpf

In preparation of the proof of Proposition
\ref{prop:OneLawConvergenceSecondMoment} we note that Proposition
\ref{prop:SecondMoment} readily implies
\begin{eqnarray}
\frac{ \bE{ T_n (\theta)^2 } }{ \bE{ T_n (\theta) }^2 } =
\frac{1}{ \bE{ T_n (\theta) } } &+& \left ( \frac{ {n-3\choose 3}
}{ {n \choose 3} } + 3 \frac{ {n-3\choose 2} }{ {n \choose 3} }
\right )
\label{eq:SecondMoment} \\
&+& \frac{ 3(n-3) }{ {n \choose 3} } \cdot \frac{
\bE{\chi_{n,123}(\theta) \chi_{n, 124} (\theta)} }
     {  \bE{\chi_{n,123}(\theta) } ^2 }
\nonumber
\end{eqnarray}
for all $n=2,3, \ldots$ as we make use of
(\ref{eq:FirstMomentExpression}).

\section{A proof of Proposition \ref{prop:OneLawConvergenceSecondMoment}}
\label{sec:ProofPropositionOneLawConvergenceSecondMoment}

Consider any scaling $P,K: \mathbb{N}_0 \rightarrow \mathbb{N}_0$
satisfying (\ref{eq:Condition1})-(\ref{eq:Condition2}). By
Proposition \ref{prop:AsymptoticEquivalence2} we have
$\lim_{n\rightarrow \infty} n^3 \beta (\theta_n) = \infty$ under
the additional condition (\ref{eq:ConditionForOne}), whence
\[
\lim_{n\rightarrow \infty} \bE{ T_n(\theta_n) } = \infty
\]
by virtue of (\ref{eq:FirstMomentExpression}).

As pointed out earlier the equivalent conditions
(\ref{eq:Condition1})-(\ref{eq:Condition2}) imply
\begin{equation}
3 K_n < P_n \label{eq:3KvsP}
\end{equation}
for all $n$ sufficiently large in $\mathbb{N}_0$. On that range
(\ref{eq:SecondMoment}) is valid with $\theta$ replaced by
$\theta_n$. Letting $n$ go to infinity in the resulting
expression, we note that
\[
\lim_{n \rightarrow \infty} \left ( \frac{ {n-3\choose 3} }{ {n
\choose 3} } + 3 \frac{ {n-3\choose 2} }{ {n \choose 3} } \right )
= 1 \quad \mbox{and} \quad \frac{ {n \choose 3} }{ 3(n-3) } \sim
\frac{n^2}{18}.
\]
It is plain that the convergence
(\ref{eq:OneLawConvergenceSecondMoment}) will hold if we show that
\begin{equation}
\lim_{n \rightarrow \infty } \frac{1}{n^2} \frac{
\bE{\chi_{n,123}(\theta_n) \chi_{n, 124} (\theta_n)} }
     { \bE{\chi_{n,123}(\theta_n) }^2 }
= 0. 
\label{eq:ToBeShown}
\end{equation}

In order to establish (\ref{eq:ToBeShown}) under the assumptions
of Proposition \ref{prop:OneLawConvergenceSecondMoment} we proceed
as follows: Recall from Lemma \ref{lem:FirstMoment} that
\begin{equation}
\bE{\chi_{n,123}(\theta_n) }^2 = \beta(\theta_n)^2 \geq \left ( 1
- q(\theta_n) \right )^6 ,
\end{equation}
and from (\ref{eq:SecondMomentCross}) observe that
\begin{eqnarray}
\lefteqn{ \frac{1}{n^2} \cdot \frac{ \bE{\chi_{n,123}(\theta_n)
\chi_{n, 124} (\theta_n)} }
           { \left ( \bE{\chi_{n,123}(\theta_n) } \right )^2 }
} & &
\nonumber \\
&=& - \frac{1}{n^2} \cdot \frac{ (1-q(\theta_n))^5 }{
\beta(\theta_n)^2 } + \frac{2}{n^2} \cdot \frac{ \left ( 1-
q(\theta_n) \right )^2 }{ \beta(\theta_n) }
\nonumber \\
& & - \frac{1}{n^2} \cdot \frac{1}{q(\theta_n)} \left ( \frac{
\beta (\theta_n) - (1-q(\theta_n))^3 }
             { \beta(\theta_n) }
\right )^2
\nonumber \\
& & + \frac{1}{n^2} \cdot \frac{ \sum_{k=0}^{K_n} c_k (\theta_n)
-q(\theta_n)^4 }{\beta(\theta_n)^2 }
\label{eq:OnelawRatioBeforeTakingLimit}
\end{eqnarray}
for all $n=3,4, \ldots $.

Let $n$ go to infinity in (\ref{eq:OnelawRatioBeforeTakingLimit}).
Using (\ref{eq:AsymptoticsForRatio}) (once with $a=5$ and $b=2$,
then with $a=2$ and $b=1$), we get
\begin{equation}
\lim_{n \rightarrow \infty} \frac{1}{n^2} \cdot \frac{
(1-q(\theta_n))^5 }{ \beta(\theta_n)^2 } = 0
\end{equation}
and
\begin{equation}
\lim_{n \rightarrow \infty} \frac{2}{n^2} \cdot \frac{ \left ( 1-
q(\theta_n) \right )^2 }{ \beta(\theta_n) } = 0 .
\end{equation}
The convergence
\begin{equation}
\lim_{n \rightarrow \infty} \frac{1}{n^2} \cdot
\frac{1}{q(\theta_n)} \left ( \frac{ \beta (\theta_n) -
(1-q(\theta_n))^3 }
             { \beta(\theta_n) }
\right )^2 = 0
\end{equation}
is immediate since
\[
\left | \frac{ \beta (\theta_n) - (1-q(\theta_n))^3 }
             { \beta(\theta_n) }
\right |^2 \leq 1, \quad n =2,3, \ldots
\]
and $\lim_{n \rightarrow \infty} q(\theta_n) = 1 $. Consequently
the proof of Proposition \ref{prop:OneLawConvergenceSecondMoment}
will be completed if we show

\begin{proposition}
{\sl For any scaling $P,K: \mathbb{N}_0 \rightarrow \mathbb{N}_0$
satisfying (\ref{eq:Condition1})-(\ref{eq:Condition2}), we have
\begin{equation}
\lim_{n \rightarrow \infty } \frac{1}{n^2} \cdot \frac{
\sum_{k=0}^K c_k (\theta_n) -q(\theta_n)^4 }{\beta(\theta_n)^2 } =
0 \label{eq:OneLawConvergenceSecondMomentReduction}
\end{equation}
under the condition (\ref{eq:ConditionForOne}). }
\label{prop:OneLawConvergenceSecondMomentReduction}
\end{proposition}

The proof of Proposition
\ref{prop:OneLawConvergenceSecondMomentReduction} will proceed in
several steps which are presented in the next three sections.

\section{The first reduction step}
\label{sec:FirstReductionStep}

We start with an easy bound.

\begin{lemma}
{\sl With positive integers $K$ and $P$ such that $2K \leq P$, 
we have
\begin{equation}
c_1(\theta)\leq 1-q(\theta) \label{eq:c_1_leq_1_q}.
\end{equation}
} 
\label{lem:c_1_leq_1_q}
\end{lemma}

\myproof Specializing (\ref{eq:SecondMomentCrossBInterpretation})
with $k=1$ we get
\begin{eqnarray*}
c_1 (\theta) &=& \bP{ |K_1(\theta) \cap K_2 (\theta) | = 1,
        \left ( K_1(\theta) \cup K_2(\theta) \right ) \cap K_i(\theta)
       = \emptyset, \ i=3,4 }
\nonumber \\
&\leq& \bP{ |K_1(\theta) \cap K_2 (\theta) | = 1 }
\nonumber \\
&\leq& \bP{ |K_1(\theta) \cap K_2 (\theta) | \geq 1 }
\end{eqnarray*}
and the conclusion is immediate as we identify
\[
\bP{ |K_1(\theta) \cap K_2 (\theta) | \geq 1 } = \bP{ K_1(\theta)
\cap K_1(\theta) \neq \emptyset } = 1 - q(\theta) .
\]
\myendpf

\begin{lemma}
{\sl With positive integers $K$ and $P$ such that $3K \leq P$, the
monotonicity property
\begin{equation}
\frac{c_1(\theta) }{c_0(\theta)} 
\geq 
\frac{c_2(\theta) }{c_1(\theta)} 
\geq \ldots \geq 
\frac{c_{K}(\theta) }{c_{K-1}(\theta)} 
\label{eq:Monotonicity}
\end{equation}
holds. 
} 
\label{lem:Monotonicity}
\end{lemma}

\myproof Fix $k=0, \ldots , K-1$. From the expression
(\ref{eq:SecondMomentCrossB}) we note that
\begin{eqnarray}
\frac{c_{k+1}(\theta)} {c_{k}(\theta)} &=& {{\begingroup {K}
\choose {k+1} \endgroup}
     {\begingroup {P-K} \choose {K-k-1} \endgroup}
     {\begingroup {P-2K+k+1} \choose {K} \endgroup}^2
\over{ \begingroup {K} \choose {k} \endgroup}
       {\begingroup {P-K} \choose {K-k} \endgroup}
       {\begingroup {P-2K+k} \choose {K} \endgroup}^2}
\nonumber \\
&=&\frac{1}{k+1} \cdot \frac{(K-k)^2}{P-3K+k+1}
\cdot\frac{P-2K+k+1}{P-3K+k+1} \label{eq:RatioExpression}
\end{eqnarray}
and by considering each factor in this last expression we readily
conclude that the ratio $\frac{c_{k+1}(\theta)}{c_{k}(\theta)}$
decreases monotonically with $k$. \myendpf

\begin{lemma}
{\sl For any scaling $P,K: \mathbb{N}_0 \rightarrow \mathbb{N}_0$
satisfying (\ref{eq:Condition1})-(\ref{eq:Condition2}), we have
\begin{equation}
\frac{c_{2}(\theta_n)}{c_{1}(\theta_n)} \leq 1 - q (\theta_n)
\label{eq:RatioVS1-q}
\end{equation}
for all $n$ sufficiently large in $\mathbb{N}_0$.}
\label{lem:RatioVS1-q}
\end{lemma}

\myproof Pick a scaling $P,K: \mathbb{N}_0 \rightarrow
\mathbb{N}_0$ satisfying
(\ref{eq:Condition1})-(\ref{eq:Condition2}) so that
(\ref{eq:3KvsP}) eventually holds. On that range replace $\theta$
by $\theta_n$ in (\ref{eq:RatioExpression}) with $k=1$ according
to this scaling, yielding
\begin{eqnarray}
\frac{c_{2}(\theta_n)}{c_{1}(\theta_n)} 
=
\frac{1}{2}\cdot\frac{(K_n-1)^2}{P_n-3K_n+2}
\cdot
\frac{P_n-2K_n+2}{P_n-3K_n+2} .
\nonumber
\end{eqnarray}
The inequality
\begin{eqnarray}
\left ( 1 - q(\theta_n) \right )^{-1}
\frac{c_{2}(\theta_n)}{c_{1}(\theta_n)}
\leq
\frac{1}{2}
\cdot
\left ( 1 - q(\theta_n) \right )^{-1}
\frac{K^2_n}{P_n-3K_n}
\cdot
\frac{P_n-2K_n}{P_n-3K_n}
\nonumber
\end{eqnarray}
readily follows.

Now let $n$ go to infinity in this inequality:
Recall the consequence
(\ref{eq:RatioConditionStrong+Consequence1}) of the assumption
(\ref{eq:Condition1})-(\ref{eq:Condition2}) and use 
the equivalence (\ref{eq:AsymptoticsEquivalence1}) to validate the
limits
\[
\lim_{n \rightarrow \infty}
\left ( 1 - q(\theta_n) \right )^{-1} \frac{K^2_n}{P_n-3K_n}
= 1
\]
and
\[
\lim_{n \rightarrow \infty}
\frac{P_n-2K_n}{P_n-3K_n}
= 1 .
\]
As a consequence,
\[
\limsup_{n \rightarrow \infty}
\left ( 1 - q(\theta_n) \right )^{-1}
\frac{c_{2}(\theta_n)}{c_{1}(\theta_n)}
\leq
\frac{1}{2}
\]
and the desired conclusion is now immediate. 
\myendpf

Combining Lemma \ref{lem:c_1_leq_1_q}, 
Lemma \ref{lem:Monotonicity} and Lemma \ref{lem:RatioVS1-q} 
will lead to the following key bounds.

\begin{lemma}
{\sl For any scaling $P,K: \mathbb{N}_0 \rightarrow \mathbb{N}_0$
satisfying (\ref{eq:Condition1})-(\ref{eq:Condition2}), we have
\begin{equation}
c_k(\theta_n)\leq\left(1-q(\theta_n)\right)^k, \quad k=1, 2,
\ldots, K_n \label{eq:c_k_leq_1-q^k}
\end{equation}
for all $n$ sufficiently large in $\mathbb{N}_0$.}
\label{lem:c_k_leq_1-q^k}
\end{lemma}

\myproof Pick a scaling $P,K: \mathbb{N}_0 \rightarrow
\mathbb{N}_0$ satisfying
(\ref{eq:Condition1})-(\ref{eq:Condition2}). For each $n=2,3,
\ldots$, we can use Lemma \ref{lem:c_1_leq_1_q} and Lemma
\ref{lem:Monotonicity} to conclude that
\begin{eqnarray}
c_k(\theta_n) &=& \prod_{\ell=1}^{k-1} \frac{c_{\ell+1}(\theta_n)}
                           {c_\ell(\theta_n)}
\cdot c_1 (\theta_n)
\nonumber \\
&\leq& \left ( \frac{c_2(\theta_n)}
             {c_1(\theta_n)} \right )^{k-1}
\cdot c_1 (\theta_n)
\nonumber \\
&\leq& \left ( \frac{c_2(\theta_n)}
             {c_1(\theta_n)} \right )^{k-1}
\cdot \left ( 1 - q(\theta_n) \right )
\end{eqnarray}
with $k=1, \ldots , K_n$. The desired conclusion is now a simple
consequence of Lemma \ref{lem:RatioVS1-q}. \myendpf

We are now in a position to take the first step towards the proof
of Proposition \ref{prop:OneLawConvergenceSecondMomentReduction}.

\begin{proposition}
{\sl For any scaling $P,K: \mathbb{N}_0 \rightarrow \mathbb{N}_0$
satisfying (\ref{eq:Condition1})-(\ref{eq:Condition2}), we have
\begin{equation}
\lim_{n \to \infty} \frac{1}{n^2} \cdot \frac{\sum_{k=5}^{K_n}
c_k(\theta_n)}
     {\beta (\theta_n)^2}
= 0 \label{eq:ReductionStep}
\end{equation}
under the condition (\ref{eq:ConditionForOne}). 
}
\label{prop:ReductionStep}
\end{proposition}

\myproof 
Pick a scaling $P,K: \mathbb{N}_0 \rightarrow \mathbb{N}_0$ 
satisfying (\ref{eq:Condition1})-(\ref{eq:Condition2}). 
The result (\ref{eq:AsymptoticsForRatio}) is trivially true 
if $K_n \leq 4$ for all $n$ sufficiently large in $\mathbb{N}_0$. 
Thus, assume from now on that $K_n \geq 5$ for infinitely many $n$ in
$\mathbb{N}_0$ -- In fact, there is now loss of generality in
assuming $K_n \geq 5$ for all $n$ sufficiently large in
$\mathbb{N}_0$. From Lemma \ref{lem:c_k_leq_1-q^k} it follows that
\begin{eqnarray}
\sum_{k=5}^{K_n} c_k(\theta_n) &\leq& \sum_{k=5}^{K_n} \left(
1-q(\theta_n) \right)^k
\nonumber \\
&\leq & \sum_{k=5}^{\infty} \left( 1-q(\theta_n) \right)^k
\nonumber \\
&=& \frac{ \left( 1-q(\theta_n) \right)^5 }{ q(\theta_n) }
\end{eqnarray}
for all $n$ sufficiently large in $\mathbb{N}_0$. Letting $n$ go
to infinity in this last inequality we readily obtain
(\ref{eq:ReductionStep}) as an immediate consequence of
Proposition \ref{prop:n^2_1_q_to_inf}, to wit
(\ref{eq:AsymptoticsForRatio}) (with $a=5$ and $b=2$). 
\myendpf

\section{The second reduction step}
\label{sec:SecondReductionStep}

It is now plain from Proposition \ref{prop:ReductionStep} that the
proof of Proposition
\ref{prop:OneLawConvergenceSecondMomentReduction} will be
completed if we show the following fact.

\begin{proposition}
{\sl For any scaling $P,K: \mathbb{N}_0 \rightarrow \mathbb{N}_0$
satisfying (\ref{eq:Condition1})-(\ref{eq:Condition2}), we have
\begin{equation}
\lim_{n \to \infty} \frac{1}{n^2} \cdot \frac{\sum_{k=0}^{4}
c_k(\theta_n) - q(\theta_n) ^ 4}
     {\beta (\theta_n)^2}
= 0 
\label{eq:SecondStep}
\end{equation}
under the condition (\ref{eq:ConditionForOne}). }
\label{prop:SecondStep}
\end{proposition}

To construct a proof of 
Proposition \ref{prop:SecondStep} 
we proceed as follows:
Fix positive integers $K$ and $P$ such that $3K \leq P$. 
By direct substitution we get
\begin{eqnarray}
\lefteqn{\sum_{k=0}^{4} c_k(\theta) - q(\theta) ^ 4 } & &
\nonumber \\
&=& \sum_{k=0}^{4} \frac{ {K \choose k}{P-K \choose K-k} }
                        { {P \choose K } } 
\left ( \frac{ {P-2K+k \choose  K} }
      { {P \choose K } } \right )^2 
- \left ( \frac{ {P-K \choose K} }{ {P \choose K } } \right )^ 4
\nonumber \\
&=& {P \choose K }^{-4} 
\left ( 
\sum_{k=0}^{4} {P \choose K } {K \choose k} {P-K \choose K-k} 
               {P-2K+k \choose  K}^2 - {P-K \choose K}^4 \right )
\nonumber \\
&=& \frac{F(\theta) }{ G(\theta) }
\label{eq:SecondStepExpression}
\end{eqnarray}
where we have set
\begin{eqnarray}
& & F(\theta)
\label{eq:F(teta)_defn} \\
&:=& 
(K!)^4 \left ( \sum_{k=0}^{4} {P \choose K } {K \choose k}
{P-K \choose K-k} {P-2K+k \choose  K}^2 - {P-K \choose K}^4 \right) 
\nonumber
\end{eqnarray}
and
\begin{equation}
G(\theta):= \left ( \frac{ P! }{ (P-K)! } \right )^4 
= \prod_{\ell=0}^{K-1}(P-\ell)^4 . 
\label{eq:G(teta)_defn}
\end{equation}

In this new notation 
Proposition \ref{prop:SecondStep} can be given a simpler, 
yet equivalent, form.

\begin{proposition}
{\sl Consider any scaling $P,K: \mathbb{N}_0 \rightarrow \mathbb{N}_0$
satisfying (\ref{eq:Condition1})-(\ref{eq:Condition2}),
The convergence (\ref{eq:SecondStep}) holds if and only if
\begin{eqnarray}
\lim_{n\to \infty}
     \frac{1}{n^2 \beta  (\theta_n) ^ 2}
     \frac{F(\theta_n)}{P_n^{4K_n}} = 0 .
\label{eq:FinalRatioToZero}
\end{eqnarray}
}
\label{prop:FinalRatioToZero}
\end{proposition}

\myproof 
Pick a scaling $P,K: \mathbb{N}_0 \rightarrow \mathbb{N}_0$ satisfying
(\ref{eq:Condition1})-(\ref{eq:Condition2}) and assume that
(\ref{eq:ConditionForOne}) holds. 
The desired equivalence is an immediate consequence
of the expression (\ref{eq:SecondStepExpression})
as we show below the equivalence
\begin{equation}
G(\theta_n) \sim P_n^{4K_n}.
\label{eq:SimpleFactForG(theta)}
\end{equation} 

By (\ref{eq:G(teta)_defn}) this last equivalence amounts to
\begin{equation}
\lim_{n \to \infty} 
\prod_{\ell=0}^{K_n-1} 
\left( \frac{P_n - \ell }{P_n} \right)^4 
= 1. 
\label{eq:P_P-l_to_1}
\end{equation}
To establish this convergence, fix $n=2,3, \ldots $ and note that
\begin{eqnarray}
\prod_{\ell=0}^{K_n-1} \left(\frac{P_n - \ell}{P_n}\right)^4 
= \left
( \prod_{\ell=0}^{K_n-1} 
\left ( 1 - \frac{\ell}{P_n} \right )
\right )^4 . 
\label{eq:Identity+A}
\end{eqnarray}
The bounds
\begin{equation}
\left( 1-\frac{K_n}{P_n} \right)^{K_n} 
\leq 
\prod_{\ell=0}^{K_n-1}
\left ( 1 - \frac{\ell}{P_n} \right ) 
\leq 1
\label{eq:boundsonProd_P_P_l}
\end{equation}
are straightforward, while simple calculus followed by a crude bounding
rgument yields
\[
1 - \left ( 1 - \frac{K_n}{P_n} \right )^{K_n} 
= \int_{ 1 - \frac{K_n}{P_n} }^1 K_n t^{K_n-1}dt 
\leq \frac{K_n^2}{P_n} .
\]
With the help of (\ref{eq:boundsonProd_P_P_l}) we now conclude that
\begin{equation}
1 -\frac{K_n^2}{P_n} 
\leq 
\prod_{\ell=0}^{K_n-1} \left( 1-\frac{\ell}{P_n} \right) 
\leq 1 .
\end{equation}
Letting $n$ go to infinity in this last expression yields
the conclusion
\begin{equation}
\lim_{n \to \infty} \prod_{\ell=0}^{K_n-1} 
\left( 1-\frac{\ell}{P_n} \right) 
= 1
\end{equation}
by virtue of (\ref{eq:Condition2}),
and this readily implies (\ref{eq:P_P-l_to_1}) via (\ref{eq:Identity+A}).
\myendpf

The following bound, which is established in Section \ref{sec:FinalStep},
proves crucial for proving the convergence
(\ref{eq:FinalRatioToZero})
under the assumptions of Proposition \ref{prop:SecondStep}.

\begin{lemma}
{\sl For any scaling $P,K: \mathbb{N}_0 \rightarrow \mathbb{N}_0$
satisfying (\ref{eq:Condition1})-(\ref{eq:Condition2}), we have
\begin{equation}
F(\theta_n) \leq K_n ^ {4} P_n ^ {4 K_n -3}
\label{eq:F_theta_bound}
\end{equation}
for all $n$ sufficiently large in $\mathbb{N}_0$.}
\label{lem:F_theta_bound}
\end{lemma}

While Lemma \ref{lem:F_theta_bound}
is established in Section \ref{sec:FinalStep},
the proof of Proposition \ref{prop:SecondStep}
can now be completed: 
Pick a scaling $P,K: \mathbb{N}_0 \rightarrow \mathbb{N}_0$ satisfying
(\ref{eq:Condition1})-(\ref{eq:Condition2}) and assume that
(\ref{eq:ConditionForOne}) holds. 
By Lemma \ref{lem:F_theta_bound} we get
\begin{equation}
 \frac{1}{n^2 \beta ^ 2 (\theta_n)} \cdot
        \frac{F(\theta_n)}{P_n^{4K_n}}
        \leq
 \frac{1}{n^2 \beta ^ 2 (\theta_n)} \cdot
        \frac{K_n^4}{P_n^{3}}
\label{eq:K^4/P^3_to_zero_pf_step1}
\end{equation}
for all $n$ sufficiently large in $\mathbb{N}_{0}$.
Invoking Proposition \ref{prop:AsymptoticEquivalence2}
we then conclude that
\begin{eqnarray}
 \frac{1}{n^2 \beta ^ 2 (\theta_n)} \cdot
        \frac{K_n^4}{P_n^{3}}
&\sim&
\frac{1}{n^2 \tau (\theta_n)^ 2 } \cdot \frac{K_n^4}{P_n^{3}}
\nonumber \\
&=& \frac{K_n^4}{n^2 P_n^3 \left( \frac{K_n^3}{P_n^2} 
+
\left(\frac{K_n^2 }{P_n}\right)^3 \right)^2 }
\nonumber \\
&\leq& \frac{K_n^4}{n^2 P_n^3 \left(\frac{K_n^3}{P_n^2}\right)^2}
\nonumber \\
&=& \left ( n^2 \frac{K_n^2}{P_n} \right )^{-1} . 
\label{eq:AsymptoticInequality}
\end{eqnarray}
The validity of (\ref{eq:FinalRatioToZero}) follows upon letting
$n$ go to infinity in (\ref{eq:K^4/P^3_to_zero_pf_step1}) and
using (\ref{eq:AsymptoticInequality}) 
together with the consequence (\ref{eq:n^2_1_q_to_infB}) of
(\ref{eq:ConditionForOne}) discussed in the proof of Proposition
\ref{prop:n^2_1_q_to_inf}. 
The proof of Proposition \ref{prop:SecondStep} is completed
with the help of Proposition \ref{prop:FinalRatioToZero}.
\myendpf

\section{Towards Lemma \ref{lem:F_theta_bound}}
\label{sec:FinalStep}

We are left with proving the key Lemma \ref{lem:F_theta_bound}.
To do so we will need to exploit the structure of $F(\theta)$:
Thus, fix positive integers $K$ and $P$ such that $3K \leq P$,
and return to (\ref{eq:F(teta)_defn}).
For each $k=0,1, \ldots , 4$, easy algebra shows that
\begin{eqnarray}
\lefteqn{ (K!)^4 {P \choose K } {K \choose k} {P-K \choose K-k}
{P-2K+k \choose  K}^2 } & &
\nonumber \\
&=& \frac{P!}{k! (P-2K+k)!} \cdot 
\left (
\frac{(K!)^2 (P-2K+k)!}{K!(K-k)! (P-3K+k)!} \right )^2 
\nonumber \\
&=& \frac{P!(P-2K+k)!}{k!} \cdot \left ( \frac{K!}{(K-k)!  (P-3K+k)!} 
\right )^2
\nonumber \\
&=& k! {K \choose k}^2 \cdot b_{K,k} (\theta)
\end{eqnarray}
with
\begin{equation}
b_{K,k} (\theta) 
:= \frac{P!(P-2K+k)!}{((P-3K+k)!)^2} .
\label{eq:b_k(theta)}
\end{equation} 
Next, it is plain that
\begin{equation}
b_K(\theta)
:= (K!)^4 {P-K \choose K}^4 
= \left ( \frac{(P-K)!}{(P-2K)!} \right)^4 .
\label{eq:b(theta)}
\end{equation}
Reporting these facts into (\ref{eq:F(teta)_defn})
we readily conclude
\begin{eqnarray}
F(\theta) 
&=& \sum_{k=0}^4 
k! {K \choose k}^2 
\cdot
\frac{P!(P-2K+k)!}{((P-3K+k)!)^2} 
- \left ( \frac{(P-K)!}{(P-2K)!} \right )^4
\nonumber \\
&=& \left ( \sum_{k=0}^4 k! {K \choose k}^2 \cdot b_{K,k}(\theta) \right )
- b_K(\theta ) .
\label{eq:F(teta)_A}
\end{eqnarray}

By direct inspection,
using (\ref{eq:Factors1}) and (\ref{eq:Factors2})
in Appendix \ref{Appendix:C},
we check that $F(\theta)$ can be written as a
polynomial in $P$ (of order $4K$), namely
\begin{equation}
F(\theta) 
= \sum_{\ell=0}^{4K} a_{4K -\ell} (K) P^\ell
= \sum_{\ell=0}^{4K} a_{\ell} (K) P^{4K-\ell}
\label{eq:F(theta)_polynom}
\end{equation}
where the coefficients 
are {\em integers} which depend on $\theta$ only through $K$. 
The first six coefficients can be evaluated explicitly.

\begin{lemma}
{\sl With positive integers $K$ and $P$ such that $3K \leq P$, we
have
\begin{equation}
a_0(K)=a_1(K)=a_2(K)=0
\end{equation}
and
\begin{equation}
a_3(K)=K^4
\label{eq:a_3}
\end{equation}
whereas
\begin{equation}
a_4(K)=-6K^6+6K^5-K^4 
\label{eq:a_4}
\end{equation}
and
\begin{eqnarray}
a_5(K) 
&=& - \frac{1}{120}K^{10} + \frac{1}{6}K^9 + \frac{199}{12}K^8
    - 34K^7 + \frac{1207}{120}K^6
\nonumber \\
& & ~ + \frac{161}{6}K^5 - \frac{209}{6}K^4 + 20K^3 -
\frac{24}{5}K^2.
\label{eq:a_5}
\end{eqnarray}
} 
\label{lem:first_six_coef}
\end{lemma}

The fact that (\ref{eq:a_5}) 
defines a polynomial expression in $K$ with rational coefficients
does not contradict the integer nature of $a_5(K)$.
In what follows we shall find it convenient to write
\begin{equation}
a_5^\star(K) 
= a_5(K) + \frac{1}{240} K^{10} .
\label{eq:a_5Star}
\end{equation}
The proof of Lemma \ref{lem:first_six_coef} is tedious and is
given in Appendix \ref{Appendix:C}. For the remaining
coefficients, we rely on the following bounds which are also derived
in Appendix \ref{Appendix:C}. 

\begin{lemma}
{\sl With positive integers $K$ and $P$ such that $3K \leq P$, we
have
\begin{equation}
| a_\ell (K) | \leq 2 \cdot (12 K^2) ^\ell, 
\quad \ell = 0, 1, \ldots, 4K. 
\label{eq:a_i_unif_bound}
\end{equation}
} 
\label{lem:a_i_unif_bound}
\end{lemma}
As expected these bounds are in agreement
witht the exact expressions obtained 
in Lemma \ref{lem:first_six_coef} for $\ell =0, 1, \ldots , 5$.

A proof of Lemma \ref{lem:F_theta_bound} can now be given:
Pick a scaling $P,K: \mathbb{N}_0 \rightarrow \mathbb{N}_0$ 
satisfying (\ref{eq:Condition1})-(\ref{eq:Condition2}) 
and replace $\theta$ by $\theta_n$ 
in (\ref{eq:F(theta)_polynom}) according to this scaling. 
As Lemma \ref{lem:first_six_coef} implies
\begin{equation}
F(\theta_n)
= K_n^4 P_n^{4K_n-3}
+ \sum_{\ell=4}^{4K_n} a_\ell(K_n) P_n^{4K_n-\ell}
\label{eq:F(theta)Reformulated}
\end{equation}
for all $n=2,3, \ldots $, 
the bound (\ref{eq:F_theta_bound}) follows if we show that
\begin{equation}
\sum_{\ell=4}^{4K_n} a_\ell(K_n) P_n^{4K_n-\ell} \leq 0
\label{eq:F_theta_bound_prf_start}
\end{equation}
for all $n$ sufficiently large in $\mathbb{N}_0$. 

To do so, apply (\ref{eq:a_i_unif_bound}) and use 
elementary arguments to get
\begin{eqnarray}
\left | \sum_{\ell=6}^{4K_n} a_\ell(K_n) P_n^{4K_n-\ell} \right |
&\leq&
\sum_{\ell=6}^{4K_n}
\left | a_\ell(K_n) \right | P_n^{4K_n-\ell}
\nonumber \\
&\leq&
\sum_{\ell=6}^{4K_n}
2 \cdot (12 K_n^2)^\ell P_n^{4K_n-\ell}
\nonumber \\
&=&
2 P_n^{4K_n}
\sum_{\ell=6}^{4K_n}
\left ( \frac{ 12 K_n^2 }{P_n} \right )^\ell
\nonumber \\
&\leq&
2 P_n^{4K_n}
\left ( \frac{ 12 K_n^2 }{P_n} \right )^6
\cdot
\sum_{\ell=0}^{\infty}
\left ( \frac{ 12 K_n^2 }{P_n} \right )^\ell
\nonumber \\
&=&
2 P_n^{4K_n}
\left ( \frac{ 12 K_n^2 }{P_n} \right )^6
\cdot
\left ( 1 - \frac{ 12 K_n^2 }{P_n} \right )^{-1}
\end{eqnarray}
for all $n$ large enough to ensure $12 K_n^2 < P_n$, say $n \geq n_1^\star$
for some finite integer $n_1^\star$; this is a simple consequence of
condition (\ref{eq:Condition1})-(\ref{eq:Condition2}).

On that range, going back to (\ref{eq:F_theta_bound_prf_start}), 
we find
\begin{eqnarray}
\lefteqn{
\sum_{\ell=4}^{4K_n} a_\ell(K_n) P_n^{4K_n-\ell}
} & &
\nonumber \\
&\leq&
a_4(K_n) P_n^{4K_n-4}
+ a_5(K_n) P_n^{4K_n-5}
+ \left | 
\sum_{\ell=6}^{4K_n} a_\ell(K_n) P_n^{4K_n-\ell}
\right | 
\nonumber \\
&\leq&
a_4(K_n) P_n^{4K_n-4}
+ a_5(K_n) P_n^{4K_n-5}
+
2 P_n^{4K_n}
\left ( \frac{ 12 K_n^2 }{P_n} \right )^6
\cdot
\left ( 1 - \frac{ 12 K_n^2 }{P_n} \right )^{-1}
\nonumber \\
&=&
P_n^{4K_n-5} \cdot L_n
\label{eq:InequalityGeom1}
\end{eqnarray}
where
\[
L_n
:= a_{4} (K_n) P_n
+ a_5(K_n)
+ 2 (12)^6 K^{10}_n \cdot \frac{K_n^2} { P_n }
\cdot
\left ( 1 - \frac{ 12 K_n^2 }{P_n} \right )^{-1} .
\]
Therefore, (\ref{eq:F_theta_bound_prf_start})
will hold for all $n$ sufficiently large in $\mathbb{N}_0$ provided
\begin{equation}
L_n \leq 0
\label{eq:InequalityToSHOW+a}
\end{equation}
for all $n$ sufficiently large in $\mathbb{N}_0$.
This last statement will be established by showing that
$L=-\infty$ where
\[
L
:= \limsup_{n \rightarrow \infty} L_n .
\]

That $L=-\infty$ can be seen as follows:
We begin with the bound
\begin{equation}
a_{4} (K_n)
= - K_n^4 (6K_n (K_n-1) +1)
\leq -K_n^4
\label{eq:InequalityC}
\end{equation}
for all $n=1,2, \ldots $. 
Next, condition (\ref{eq:Condition1})-(\ref{eq:Condition2})
implies
\begin{equation}
\lim_{n \rightarrow \infty}
\frac{K_n^2} { P_n }
\cdot
\left ( 1 - \frac{ 12 K_n^2 }{P_n} \right )^{-1}
= 0 ,
\label{eq:Limit=Zero}
\end{equation}
whence there exists some finite integer $n_2^{\star}$ such that
\begin{equation}
2 (12)^6 
\frac{K_n^2} { P_n }
\cdot
\left ( 1 - \frac{ 12 K_n^2 }{P_n} \right )^{-1}
\leq \frac{1}{240} ,
\quad n \geq n_2^{\star} .
\label{eq:InequalityD}
\end{equation}

Now, set $n^\star = \max \left ( n_1^\star, n_2^\star \right )$,
and recall the definition (\ref{eq:a_5Star}).
On the range $n \geq n^\star$, 
both inequalities (\ref{eq:InequalityGeom1}) and
(\ref{eq:InequalityD}) hold, and we obtain 
\begin{eqnarray}
\lefteqn{a_{4} (K_n) P_n + a_5(K_n)
+ 2 (12)^6 K^{10}_n \cdot \frac{K_n^2} { P_n }
\cdot
\left ( 1 - \frac{ 12 K_n^2 }{P_n} \right )^{-1}  } & &
\nonumber \\
&=& a_{4} (K_n) P_n + a_5^\star (K_n)
+ \left ( 
- \frac{1}{240}
+ 2 (12)^6 \cdot \frac{K_n^2} { P_n }
\cdot
\left ( 1 - \frac{ 12 K_n^2 }{P_n} \right )^{-1}
\right )
K_n^{10}
\nonumber \\
&\leq& -K_n^4 P_n
+ a_5^\star (K_n)
\label{eq:InequalityE}
\end{eqnarray}
upon making use of (\ref{eq:InequalityC}). 
To conclude, set
\begin{equation}
L^\star := 
\limsup_{n \rightarrow \infty} 
\left ( a_5^\star (K_n) \right )
\label{eq:InequalityF}
\end{equation}
and note that $L^\star$ is necessarily an element of $[-\infty, \infty )$,
i.e., it is never the case that $L^\star = \infty$.
This follows easily from the fact that the mapping
$\mathbb{R}_+ \rightarrow \mathbb{R}_+: 
x \rightarrow a_5^\star (x) $
is a polynomial of degree $10$ whose leading coefficient 
($-\frac{1}{240}$) is negative.
As we recall (\ref{eq:RatioConditionStrong+Consequence2})
under (\ref{eq:Condition1})-(\ref{eq:Condition2}),
it is now plain from (\ref{eq:InequalityE})
that $L = -\infty$ by standard properties of the lim sup operation.
\myendpf

Careful inspection of the proof 
of Proposition \ref{prop:SecondStep}
given at the end of Section \ref{sec:SecondReductionStep}
shows that the inequality (\ref{eq:F_theta_bound}) of
Lemma \ref{lem:F_theta_bound} could be replaced 
without prejudice by the following weaker statement:
For any scaling $P,K: \mathbb{N}_0 \rightarrow \mathbb{N}_0$
satisfying (\ref{eq:Condition1})-(\ref{eq:Condition2}), there
exists some positive constant $C$ such that
\begin{equation}
F(\theta_n) \leq C K_n ^ {4} P_n ^ {4 K_n -3}
\label{eq:F_theta_bound+Weaker}
\end{equation}
for all $n$ sufficiently large in $\mathbb{N}_0$.

Now, from only the knowledge of the first four coefficients
in Lemma \ref{lem:first_six_coef} we can already conclude that
\begin{equation}
\lim_{P \rightarrow \infty}
\frac{F(K,P)}{K^4 P^{4K-3} } = 1
\end{equation}
for {\em each} $K=1,2, \ldots $, so that for each $\varepsilon > 0$
there exists a finite integer $P^\star ( \varepsilon, K) $ such that
\begin{equation}
F(K,P) \leq \left ( 1 + \varepsilon \right ) K^4 P^{4K-3},
\quad P \geq P^\star ( \varepsilon, K) 
\end{equation}
Unfortunately, the threshold 
$P^\star ( \varepsilon, K) $ is not known to be 
uniform with respect to $K$, and the approach 
does {\em not} necessarily imply
(\ref{eq:F_theta_bound+Weaker}) (with $C = 1 + \varepsilon$)
{\em unless} the sequence $K: \mathbb{N}_0 \rightarrow \mathbb{N}_0$
is bounded. This technical difficulty
is at the root of why more information on
the coefficients $a_4(K)$ and $a_5(K)$ 
(as provided in Lemma \ref{lem:first_six_coef})
is needed, and paves the way for the subsequent arguments behind
Lemma \ref{lem:F_theta_bound}.

\appendix
\setcounter{equation}{0}
\renewcommand{\theequation}{\thesection.\arabic{equation}}

\section{Establishing (\ref{eq:AsymptoticsEquivalence2Reduced})}
\label{Appendix:A}

With positive integers $K, P$  such that $3K \leq P$, we note that
\begin{eqnarray}
\frac{r(\theta)}{q(\theta)^2}
&=& 
\left( \frac { ( P - 2 K ) ! } { ( P - K ) ! } \right ) ^ 2
\cdot 
\frac { ( P - 2 K ) ! } { ( P - 3 K ) !  } 
\cdot \frac { P !  } { ( P - K ) ! }
\nonumber \\
&=& \prod_{\ell=0}^{K-1} 
\left ( \frac{P-2K-\ell}{P-K-\ell} \right) 
\cdot 
\prod_{\ell=0}^{K-1} \left ( \frac{P-\ell}{P-K-\ell} \right )
\nonumber \\
&=& \prod_{\ell=0}^{K-1} 
\left( 1 - \left(\frac{K}{P-K-\ell}\right)^2 \right) ,
\end{eqnarray}
and elementary bounding arguments yield
\[
\left(1-\left(\frac{K}{P-2K}\right)^2\right)^K 
\leq \frac{ r ( \theta ) } { q ( \theta ) ^ 2 } 
\leq
\left(1-\left(\frac{K}{P-K}\right)^2\right)^K .
\]

Pick a scaling $P,K: \mathbb{N}_0 \rightarrow \mathbb{N}_0$
satisfying the equivalent conditions
(\ref{eq:Condition1})-(\ref{eq:Condition2}) and consider $n$
sufficiently large in $\mathbb{N}_0$ so that (\ref{eq:Condition0})
holds with $c=3$. On that range, as we replace $\theta $ by
$\theta_n $ in the last chain of inequalities according
to this scaling, we get
\[
1 - \left( 1-\left( \frac{K_n}{P_n-K_n} \right )^2 \right)^{K_n}
\leq 1 - \frac{ r ( \theta_n ) } { q ( \theta_n )^ 2} \leq 1 -
\left( 1-\left( \frac{K_n}{P_n-2K_n} \right)^2 \right)^{K_n}.
\]
A standard sandwich argument will imply the desired
equivalence (\ref{eq:AsymptoticsEquivalence2Reduced})
if we show that
\begin{equation}
1 - \left( 1-\left( \frac{K_n}{P_n-cK_n} \right)^2 \right)^{K_n}
\sim \ \frac{K_n^3}{P^2_n}, \quad c=1,2 . 
\label{eq:AsymptoticEquivToBeShown}
\end{equation}

To establish (\ref{eq:AsymptoticEquivToBeShown}) 
we proceed as follows: Fix $c=1,2$ 
and on the appropriate range we note that
\begin{eqnarray}
\lefteqn{ 1 - \left( 1-\left( \frac{K_n}{P_n-cK_n} \right)^2
\right)^{K_n} } & &
\nonumber \\
&=& \int_{1- \left ( \frac{K_n}{P_n-cK_n} \right)^2 }^1 K_n
t^{K_n-1} dt
\nonumber \\
&=& K_n \left( \frac{K_n}{P_n-cK_n} \right)^2 \int_{0}^1 \left ( 1
- \left( \frac{K_n}{P_n-cK_n} \right)^2 \tau \right )^{K_n-1} d\tau 
\label{eq:BasicEquality}
\end{eqnarray}
after performing the simple change of variables $t = 1 - \left(
\frac{K_n}{P_n-cK_n} \right)^2 \tau$.

Next we invoke (\ref{eq:RatioConditionStrong+Consequence1}) to
find
\begin{equation}\label{eq:K/P-K_K/P}
\left ( \frac{K_n}{P_n-cK_n} \right)^2 = \left (
\frac{K_n}{P_n}\left(1+o(1)\right) \right )^2
=\frac{{K_n}^2}{{P_n}^2}(1+o(1)) ,
\end{equation}
so that
\begin{equation}
K_n\left(\frac{K_n}{P_n-cK_n}\right)^2 \sim \frac{K^3_n }{P^2_n} .
\label{eq:BasicEquality2}
\end{equation}
It is now plain from (\ref{eq:BasicEquality}) and
(\ref{eq:BasicEquality2}) that (\ref{eq:AsymptoticEquivToBeShown}) 
holds provided
\begin{equation}
\lim_{n \rightarrow \infty} \int_{0}^1 \left ( 1 - \left(
\frac{K_n}{P_n-cK_n} \right)^2 \tau \right )^{K_n-1} d\tau = 1 .
\end{equation}
This is a consequence of the Bounded Convergence Theorem since
\[
\lim_{n \rightarrow \infty} \left ( 1 - \left(
\frac{K_n}{P_n-cK_n} \right)^2 \tau \right )^{K_n-1} = 1, \quad 0
\leq \tau \leq 1
\]
upon noting by elementary convergence results that
\[
\lim_{n \to \infty} K_n\left(\frac{K_n}{P_n-cK_n}\right)^2 \tau =
\lim_{n \to \infty} \left ( \frac{K^2_n }{P_n} \right ) \left (
\frac{K_n}{P_n} \right ) \tau = 0
\]
across the range as a direct consequence of (\ref{eq:Condition2})
and (\ref{eq:RatioConditionStrong+Consequence1}).  
\myendpf

\setcounter{equation}{0}
\renewcommand{\theequation}{\thesection.\arabic{equation}}
\section{Evaluating
         (\ref{eq:SecondMomentCross})--(\ref{eq:SecondMomentCrossB})}
\label{Appendix:B}

For notational convenience, we define
\[
K_{ij}:=[K_i(\theta) \cap K_j(\theta) \neq \emptyset].
\]
for distinct $i,j = 1, 2, \ldots n$. Moreover, for any non-empty
subset $S$ of $\{ 1, \ldots , P \}$, we write
\[
K_{Si}:=[S \cap K_i(\theta) \neq \emptyset], \quad i=1, \ldots ,
n.
\]

In what follows we make repeated use of the decomposition
(\ref{eq:Set_Formula}). Beginning with the observation
\begin{eqnarray}
\lefteqn{\bE{\chi_{n,123}(\theta) \chi_{n, 124} (\theta)} } & &
\nonumber \\
&=& \bP{K_{12}, K_{13}, K_{23}, K_{14}, K_{24}}
\nonumber \\
&=& \bP{  K_{13}, K_{23}, K_{14}, K_{24}}
     - \bP{K_{12}^{c}, K_{13}, K_{23}, K_{14}, K_{24}}.
\label{eq:Start_Prop}
\end{eqnarray}
we shall compute each term in turn.

To compute the second term in (\ref{eq:Start_Prop}) we condition
on the sets $K_1$ and $K_2$ such that $ K_1 \cap K_2 = \emptyset$.
Thus,
\begin{eqnarray}\nonumber
\lefteqn{ \bP {K_{12}^{c}, K_{13}, K_{23}, K_{14}, K_{24}} } & &
\nonumber \\
&=& \sum_{|S|=|T|=K,  S\cap T = \emptyset} \bP{K_{1}=S,K_{2}=T,
K_{S3},K_{T3},K_{S4},K_{T4}}
\nonumber \\
&=& \sum_{|S|=|T|=K, S\cap T = \emptyset} \bP{K_{1}=S,K_{2}=T}
\bP{K_{S3}, K_{T3}, K_{S4}, K_{T4}}
\nonumber \\
&=& \sum_{|S|=|T|=K,  S\cap T = \emptyset} {P \choose K}^{-2}
\bP{K_{S3}, K_{T3}}\cdot\bP{K_{S4}, K_{T4}}
\nonumber \\
&=& \sum_{|S|=|T|=K,  S\cap T = \emptyset} {P \choose K}^{-2}
\left ( \bP{K_{S3}, K_{T3}} \right )^2
\nonumber \\
&=& {P \choose K}^{-2} \sum_{|S|=|T|=K,  S\cap T = \emptyset}
\left( \bP{K_{S3}} - \bP{K_{T3}^{c}} + \bP{K_{S3}^{c}, K_{T3}^c}
\right)^2
\nonumber \\
&=& {P \choose K}^{-2} \sum_{|S|=|T|=K,  S\cap T = \emptyset}
\left( 1 - \bP{K_{S3}^{c}} - \bP{K_{T3}^{c}} + \bP{K_{S3}^{c},
K_{T3}^c} \right)^2
\nonumber \\
&=& {P \choose K}^{-2} \sum_{|S|=|T|=K,  S\cap T = \emptyset}
\left( 1- 2 q(\theta) + r(\theta) \right )^2
\nonumber \\
&=& {P \choose K}^{-2} {P \choose K}{{P-K}\choose K} \left( 1- 2
q(\theta) + r(\theta) \right )^2
\nonumber \\
&=& q(\theta) \left( 1- 2 q(\theta) + r(\theta) \right )^2
\label{eq:subs1}
\end{eqnarray}
as we note from (\ref{eq:Probab_key_ring_does_not_intersect_S})
that $\bP{K_{S3}^{c}} = \bP{K_{T3}^{c}} = q(\theta)$ for $S$ and
$T$ in ${\cal P}_K$ with $\bP{K_{S3}^{c}, K_{T3}^c} = r(\theta)$
whenever $S \cap T = \emptyset$.

We now turn to the first term in (\ref{eq:Start_Prop}). Again,
upon making repeated use of (\ref{eq:Set_Formula}) we find
\begin{eqnarray}
\lefteqn{\bP{K_{13}, K_{23}, K_{14}, K_{24}}} & &
\nonumber \\
&=& \bP{ {K_{23}, K_{14}, K_{24}}} - \bP{ K_{13}^{c}, K_{23},
K_{14}, K_{24}}
\nonumber \\
&=& \bP{ K_{14}, K_{24}} - \bP {K_{23}^{c}, K_{14}, K_{24}} - \bP{
K_{13}^{c}, K_{14}, K_{24}} + \bP{ K_{13}^{c}, K_{23}^{c}, K_{14},
K_{24}}
\nonumber \\
&=& (1-q(\theta))^2 - 2 \bP {K_{23}^{c}, K_{14}, K_{24}}
 + \bP{  K_{13}^{c}, K_{23}^{c}, K_{24}}
 - \bP{  K_{13}^{c}, K_{23}^{c}, K_{14}^{c}, K_{24} }
\nonumber \\
&=& (1-q(\theta))^2 - 2 \bP {K_{23}^{c}, K_{14}, K_{24}}
    + \bP{ K_{13}^{c}, K_{23}^{c}, K_{24}}
\nonumber \\ \label{eq:subs2} & & ~ - \bP{ K_{13}^{c}, K_{23}^{c},
K_{14}^{c} }
      + \bP{ K_{13}^{c}, K_{23}^{c}, K_{14}^{c}, K_{24}^{c}}
\end{eqnarray}
as we note that $\bP {K_{23}^{c}, K_{14}, K_{24}} 
= \bP{ K_{13}^{c}, K_{14}, K_{24}}$. Next, we find
\begin{eqnarray}
\bP {K_{23}^{c}, K_{14}, K_{24}} &=& \sum_{|S|=K} \bP{K_4 = S,
K_{23}^{c}, K_{S1}, K_{S2} }
\nonumber \\
&=& \sum_{|S|=K} \bP{K_{4}=S} \bP{K_{23}^{c}, K_{S1}, K_{S2}}
\nonumber \\
&=& \sum_{|S|=K} {P \choose K}^{-1} \bP{K_{S1}} \cdot
\bP{K_{23}^{c}, K_{S2}}
\nonumber \\
&=& \sum_{|S|=K} {P \choose K}^{-1}(1-q(\theta))
    \cdot q(\theta)(1-q(\theta))
\label{eq:use_q_1-q} \\
&=& q(\theta)(1-q(\theta))^2 \label{eq:subs3}
\end{eqnarray}
upon using (\ref{eq:cor_q_1-q}) in (\ref{eq:use_q_1-q}).

In a similar manner, we obtain
\begin{eqnarray}
\bP {K_{13}^{c}, K_{23}^{c}, K_{24}} &=& \sum_{|S|=K} \bP{K_{2}=S,
K_{13}^{c}, K_{S3}^{c}, K_{S4} }
\nonumber \\
&=& \sum_{|S|=K} \bP{K_{2}=S} \bP{K_{13}^{c}, K_{S3}^{c}, K_{S4}}
\nonumber \\
&=& \sum_{|S|=K} {P \choose K}^{-1} \bP{K_{S4}} \cdot
\bP{K_{13}^{c}, K_{S3}^{c}}
\nonumber \\
&=&\sum_{|S|=K} {P \choose K}^{-1} (1-q(\theta))\cdot q(\theta)^2
\label{eq:use_q^2} \\
&=& q(\theta)^2 (1-q(\theta)) \label{eq:subs4}
\end{eqnarray}
where (\ref{eq:use_q^2}) follows from (\ref{eq:cor_q^2}).

We also get
\begin{eqnarray}
\bP {K_{13}^{c}, K_{23}^{c}, K_{14}^{c}} &=& \sum_{|S|=K}
\bP{K_{1}=S, K_{S3}^{c}, K_{23}^{c}, K_{S4}^{c} }
\nonumber \\
&=& \sum_{|S|=K} \bP{K_{1}=S} \bP{K_{S3}^{c}, K_{23}^{c},
K_{S4}^{c}}
\nonumber \\
&=& \sum_{|S|=K} {P \choose K}^{-1} \bP{K_{S4}^{c}} \cdot
\bP{K_{S3}^{c}, K_{23}^{c}}
\nonumber \\
&=& \sum_{|S|=K} {P \choose K}^{-1} q(\theta) \cdot q(\theta)^2
\nonumber \\
&=& q(\theta)^3. \label{eq:subs5}
\end{eqnarray}

Finally consider the term $\bP{ K_{13}^{c}, K_{23}^{c},
K_{14}^{c}, K_{24}^{c}}$: By conditioning on the cardinality of
the intersection $K_1 \cap K_2$, we obtain
\begin{eqnarray}
\lefteqn{ \bP{ K_{13}^{c}, K_{23}^{c}, K_{14}^{c}, K_{24}^{c}}} &
&
\nonumber \\
&=& \sum_{|S|=|T|=K} \bP{K_{1}=S,K_{2}=T,
K_{S3}^{c},K_{T3}^{c},K_{S4}^{c},K_{T4}^{c}}
\nonumber \\
&=& \sum_{|S|=K} \sum_{k=0}^{K} \sum_{|T|=K, |T \cap S|=k}
\bP{K_{1}=S,K_{2}=T, K_{S3}^{c},K_{T3}^{c},K_{S4}^{c},K_{T4}^{c}}
\nonumber \\
&=& \sum_{k=0}^{K} c_k (\theta)
\end{eqnarray}
where for each $k=0,1, \ldots , K$, we have set
\begin{eqnarray}
c_k(\theta) &:=& \sum_{|S|=|T|=K, |T \cap S|=k}
\bP{K_{1}=S,K_{2}=T, K_{S3}^{c},K_{T3}^{c},K_{S4}^{c},K_{T4}^{c}}
\label{eq:subs5to6} \\
&=& \sum_{|S|=|T|=K, |T \cap S|=k} \bP{K_{1}=S} \bP{K_2=T}
\bP{K_{S3}^{c}, K_{T3}^{c}} \cdot \bP{ K_{S4}^{c}, K_{T4}^{c}}
\nonumber \\
&=& \sum_{|S|=K} \bP{K_{1}=S} \sum_{|T|=K, |T \cap S|=k}
\bP{K_2=T} \bP{K_{S3}^{c}, K_{T3}^{c}} \cdot \bP{ K_{S4}^{c},
K_{T4}^{c}}
\nonumber \\
&=& \sum_{|S|=K} {P \choose K}^{-1} \sum_{|T|=K, |T \cap S|=k} {P
\choose K}^{-1} \bP{K_{S3}^{c}, K_{T3}^{c}} \cdot \bP{ K_{S4}^{c},
K_{T4}^{c}}
\nonumber \\
&=& \sum_{|S|=K} {P \choose K}^{-1} \sum_{|T|=K, |T \cap S|=k} {P
\choose K}^{-1} \left ( \frac{ {P-2K+k \choose K} }
     { {P \choose K} }
\right )^2
\nonumber \\
&=& \sum_{|S|=K} {P \choose K}^{-1} \cdot {K \choose k}{P-K
\choose K-k} \cdot {P \choose K}^{-1} \left ( \frac{ {P-2K+k
\choose K} }
     { {P \choose K} }
\right )^2 \label{eq:subs6}
\end{eqnarray}
and the expression follows (\ref{eq:SecondMomentCrossB}).

Substituting (\ref{eq:subs1}) and (\ref{eq:subs2}) (with the help
of (\ref{eq:subs3}), (\ref{eq:subs4}), (\ref{eq:subs5}) and
(\ref{eq:subs6})) into (\ref{eq:Start_Prop}), we find
\begin{eqnarray}
\lefteqn{ \bE{\chi_{n,123}(\theta) \chi_{n, 124} (\theta)} } & &
\nonumber \\
&=& (1-q(\theta))^2 - 2 q(\theta) (1-q(\theta))^2 + q(\theta)^2
(1-q(\theta)) - q(\theta)^3
\nonumber \\
& & - q(\theta) \left( 1- 2 q(\theta) + r(\theta) \right )^2
    + \sum_{k=0}^K c_k (\theta)
\label{eq:SecondMomentCrossA}
\end{eqnarray}
where we have used the notation (\ref{eq:SecondMomentCrossB}).

As we seek to simplify this last expression, we  note that
\begin{eqnarray}
\lefteqn{ (1-q(\theta))^2 - 2 q(\theta) (1-q(\theta))^2 +
q(\theta)^2 (1-q(\theta)) - q(\theta)^3 } & &
\nonumber \\
&=& (1-q(\theta))^2 \left ( 1 - 2 q(\theta) \right )
    + q(\theta)^2 (1-q(\theta)) - q(\theta)^3
\nonumber \\
&=& (1-q(\theta))^2 \left ( 1 - 2 q(\theta)  + q(\theta)^2 \right)
    - q(\theta)^2 (1-q(\theta))^2
\nonumber \\
& & ~ + q(\theta)^2 (1-q(\theta)) - q(\theta)^3
\nonumber \\
&=& (1-q(\theta))^4
    + q(\theta)^2 \left ( (1-q(\theta)) - (1-q(\theta))^2 \right )
     - q(\theta)^3
\nonumber \\
&=& (1-q(\theta))^4
    + q(\theta)^2 (1- q(\theta)) \left ( 1 - (1-q(\theta)) \right )
     - q(\theta)^3
\nonumber \\
&=& (1-q(\theta))^4
    + q(\theta)^3 (1- q(\theta))
     - q(\theta)^3
\nonumber \\
&=& (1-q(\theta))^4 - q(\theta)^4 . \label{eq:Calculations1}
\end{eqnarray}
Next, we observe that
\begin{eqnarray}
\lefteqn{ q(\theta) \left( 1- 2 q(\theta) + r(\theta) \right )^2 }
& &
\nonumber \\
&=& q(\theta) \left( 1- 2 q(\theta) + q(\theta)^2 - q(\theta)^2
                     + r(\theta) \right )^2
\nonumber \\
&=& q(\theta) \left( \left ( 1- q(\theta) \right )^2
                 - \left ( q(\theta)^2 - r(\theta) \right ) \right )^2
\nonumber \\
&=& q(\theta) \left( \left ( 1- q(\theta) \right )^4 - 2 \left (
1- q(\theta) \right )^2 \left ( q(\theta)^2 - r(\theta) \right ) +
\left ( q(\theta)^2 - r(\theta) \right )^2 \right )
\nonumber \\
&=& q(\theta) \left ( 1- q(\theta) \right )^4 - 2 q(\theta) \left
( 1- q(\theta) \right )^2
              \left ( q(\theta)^2 - r(\theta) \right )
\nonumber \\
& & ~ + q(\theta) \left ( q(\theta)^2 - r(\theta) \right )^2 .
\label{eq:Calculations2}
\end{eqnarray}

Subtracting (\ref{eq:Calculations2}) from (\ref{eq:Calculations1})
gives
\begin{eqnarray}
\lefteqn{ (1-q(\theta))^4 - q(\theta)^4 - q(\theta) \left( 1- 2
q(\theta) + r(\theta) \right )^2 } & &
\nonumber \\
&=& (1-q(\theta))^4 - q(\theta)^4 - q(\theta) \left ( 1- q(\theta)
\right )^4 + 2 q(\theta) \left ( 1- q(\theta) \right )^2 \left (
q(\theta)^2 - r(\theta) \right )
\nonumber \\
& & ~ - q(\theta) \left ( q(\theta)^2 - r(\theta) \right )^2
\nonumber \\
&=& (1-q(\theta))^4 \left ( 1 - q(\theta) \right ) - q(\theta)^4 +
2 q(\theta) \left ( 1- q(\theta) \right )^2
              \left ( q(\theta)^2 - r(\theta) \right )
\nonumber \\
& & ~ - q(\theta) \left ( q(\theta)^2 - r(\theta) \right )^2
\nonumber \\
&=& (1-q(\theta))^5 - q(\theta)^4 + 2 q(\theta) \left ( 1-
q(\theta) \right )^2
              \left ( q(\theta)^2 - r(\theta) \right )
\nonumber \\
& & ~ - q(\theta) \left ( q(\theta)^2 - r(\theta) \right )^2 .
\end{eqnarray}
Reporting the outcome of this last calculation into
(\ref{eq:SecondMomentCrossA}) we then get
\begin{eqnarray}
\lefteqn{ \bE{\chi_{n,123}(\theta) \chi_{n, 124} (\theta)} } & &
\nonumber \\
&=& (1-q(\theta))^5 + 2 q(\theta) \left ( 1- q(\theta) \right )^2
              \left ( q(\theta)^2 - r(\theta) \right )
\nonumber \\
& & ~ - q(\theta) \left ( q(\theta)^2 - r(\theta) \right )^2
      + \sum_{k=0}^K c_k (\theta) - q(\theta)^4,
\end{eqnarray}
and the conclusion (\ref{eq:SecondMomentCross}) follows as we make
use of the expression (\ref{eq:beta(tetha)}) for $\beta(\theta)$.

\setcounter{equation}{0}
\renewcommand{\theequation}{\thesection.\arabic{equation}}
\section{Proofs of Lemma \ref{lem:first_six_coef} and 
Lemma \ref{lem:a_i_unif_bound}}
\label{Appendix:C}

The proofs of both
Lemma \ref{lem:first_six_coef} and 
Lemma \ref{lem:a_i_unif_bound} will make use of the
following observations:
Pick positive integers $K $ and $P$ such that $K \geq 4$ and $3K \leq P$, 
and recall the expressions (\ref{eq:b_k(theta)}) and 
(\ref{eq:b(theta)}) appearing in (\ref{eq:F(teta)_A}).

Fix $k=0,1, \ldots , 4$. The product
\begin{equation}
b_{K,k}(\theta)
=
\prod_{i=1}^{K} (P-3K+k+i)
\cdot 
\prod_{j=1}^{3K-k} (P-3K+k+j)
\label{eq:Factors1}
\end{equation}
has $K + (3K-k) = 4K-k$ factors, hence
defines a polynomial expression in $P$
with leading term $P^{4K-k}$, say
\begin{equation}
b_{K,k}(P) = \sum_{\ell=0}^{4K-k} \beta_{k,\ell}(K) P^\ell
\end{equation}
for some integer coefficients 
$\beta_{k,0}(K) , \ldots , \beta_{k,4K-k}(K)$
with $\beta_{k,4K-k}(K) = 1$.  
On the other hand, the expression
\begin{equation}
b_K(\theta) 
=
\left ( \prod_{i=1}^K (P-2K+i) \right )^4 
= 
\left ( \prod_{i=K}^{2K-1} (P - i ) \right )^4 
\label{eq:Factors2}
\end{equation}
is a product of $4K$ factors with leading term $P^{4K}$,
and we can write it as a polynomial in $P$, namely
\begin{equation}
b_K(P) = \sum_{\ell=0}^{4K} \beta_{\ell}(K) P^\ell
\end{equation}
for some integer coefficients 
$\beta_{0}(K) , \ldots , \beta_{4K}(K)$ with $\beta_{4K}(K) = 1$.

Direct substitution followed by elementary manipulations gives
\begin{eqnarray}
\lefteqn{
\sum_{k=0}^4 k! {K \choose k}^2 \cdot b_{K,k}(\theta) } & &
\nonumber \\
&=& \sum_{k=0}^4 k! {K \choose k}^2 \cdot 
\left ( \sum_{\ell=0}^{4K-k} \beta_{k,\ell}(K) P^\ell \right )
\nonumber \\
&=& \sum_{k=0}^4 k! {K \choose k}^2 \cdot 
\left ( \sum_{\ell=0}^{4K-5} \beta_{k,\ell}(K) P^\ell 
+
\sum_{\ell=4K-4}^{4K-k} \beta_{k,\ell}(K) P^\ell 
\right )
\nonumber \\
&=&
\sum_{\ell=0}^{4K} 
\left (
\sum_{k=0}^{\min(4K-\ell,4)} k! {K \choose k}^2 
\beta_{k,\ell}(K) \right ) P^\ell ,
\nonumber
\end{eqnarray}
and it is then plain that
\begin{eqnarray}
F(\theta)
&=& 
\sum_{k=0}^4 k! {K \choose k}^2 \cdot b_{K,k}(\theta) - b(\theta)
\nonumber \\
&=& 
\sum_{\ell=0}^{4K} 
\left (
\sum_{k=0}^{\min(4K-\ell,4)} k! {K \choose k}^2 
\beta_{k,\ell}(K) \right ) P^\ell
- \sum_{\ell=0}^{4K} \beta_{\ell}(K) P^\ell .
\nonumber 
\end{eqnarray}
Finally, upon comparing with (\ref{eq:F(theta)_polynom}) we get
the relations
\begin{equation}
a_{\ell} (K)
= \left (
\sum_{k=0}^{\min(\ell,4)} k! {K \choose k}^2 \beta_{k,4K-\ell}(K) \right ) 
- \beta_{4K-\ell}(K) 
\label{eq:CoefficientFormula+Equiv}
\end{equation}
for all $\ell=0, \ldots , 4K$.

\subsection{A proof of Lemma \ref{lem:a_i_unif_bound}}

We begin with some simple observations: For some positive integer
$M$, consider the mapping $R: \mathbb{R} \rightarrow \mathbb{R}$
given by
\[
R(x) = \prod_{m=1}^M \left ( x - r_m \right ),
\quad  x \in \mathbb{R}
\]
with scalars $r_1, \ldots , r_M$, not necessarily distinct.
Obviously, $R: \mathbb{R} \rightarrow \mathbb{R}$ is a polynomial
(in the variable $x$)
of degree $M$ with all its roots located at $r_1, \ldots , r_M$.
It can be written in the form
\begin{equation}
R(x) = \sum_{m=0}^M \rho_{M-m} x^m,
\quad  x \in \mathbb{R}
\end{equation}
for some coefficients $\rho_0 , \ldots , \rho_M$ with $\rho_0 =1$;
these coefficients are uniquely 
determined by the roots $r_1, \ldots , r_M$.
In fact, for each $m=0,1, \ldots , M$, 
the coefficient $\rho_m$ of $x^{M-m}$ is given by
\begin{equation}
\rho_{m} = (-1)^{m} {\sum}_{(k_1, \ldots , k_m) \in {\cal M}_m}
 r_{k_1} \ldots r_{k_m}
\label{eq:FromRootsToCoefficients}
\end{equation}
where ${\cal M}_m$ denotes the collection of all 
unordered $m$-uples drawn without repetition
from the set of indices $1, \ldots , M$.
Obviously $| {\cal M}_m | = { M \choose m } $ and the bounds
\begin{equation}
\left | \rho_m \right |
\leq 
{ M \choose m } \cdot \left ( r^\star \right )^m
\leq 
\left ( M r^\star \right )^m
\label{eq:GenericBoundsCoefficients1}
\end{equation}
hold with $r^\star$ given by
\begin{equation}
r^\star
:= \max \left ( |r_m|, \ m=1, \ldots , M \right ) .
\label{eq:GenericBoundsCoefficients2}
\end{equation}

Now we turn to the proof of Lemma \ref{lem:a_i_unif_bound}:
Pick positive integers $K$ and $P$ such that $K \geq 4$ and
$3K \leq P$, and fix $\ell=4,5, \ldots , 4K$ -- We shall give
a proof only in that range for simplicity of exposition;
after all the desired bounds are already implied by the exact
expression for $a_0(K), \ldots , a_3(K)$ given as part of
Lemma \ref{lem:first_six_coef}.
On the range $\ell =4, \ldots , 4K$,
the bound (\ref{eq:CoefficientFormula+Equiv}) already implies
\begin{equation}
\left | a_{\ell} (K) \right |
\leq
\left (
\sum_{k=0}^{4} k! {K \choose k}^2 
\left | \beta_{k,4K - \ell}(K) \right | \right ) 
+ \left | \beta_{4K-\ell}(K)  \right | .
\label{eq:InequalityForAbsolutevalue}
\end{equation}
For each $k=0, 1, \ldots , 4$, we note that
\begin{equation}
k! {K \choose k}^2
\leq 
k! \left ( \frac{K^k}{k!} \right )^2
= \frac{K^{2k}}{k!}.
\label{eq:Bound+X}
\end{equation}
We then apply the bound
(\ref{eq:GenericBoundsCoefficients1})-(\ref{eq:GenericBoundsCoefficients2}) 
to the polynomal $b_{K,k}$:
From (\ref{eq:Factors1}) we get the values $M = 4K-k$
and $r^\star = 3K -(k+1)$, Also,
we note that $\beta_{k,4K - \ell}(K)$ is 
the coefficient of $P^{4K-\ell} $ (thus of $P^{4K - k -(\ell -k)} $)
in the polynomial $b_{K,k}(P)$ of order $4K-k$.
Therefore,
applying the bound
(\ref{eq:GenericBoundsCoefficients1})-(\ref{eq:GenericBoundsCoefficients2}) 
we find
\begin{eqnarray}
\left | \beta_{k,4K - \ell}(K) \right |
&\leq &
\left ( (4K-k  ) \cdot ( 3K-(k+1) ) \right )^{\ell-k}
\nonumber \\
&\leq& \left ( 12 K^2 \right )^{\ell - k}.
\label{eq:Bound+Y}
\end{eqnarray}
In a similar way, we apply the bound
(\ref{eq:GenericBoundsCoefficients1})-(\ref{eq:GenericBoundsCoefficients2}) 
to the polynomial $b_K$. 
This time, (\ref{eq:Factors2}) gives $M = 4K$
and $r^\star = 2K -1$, and we conclude that
\begin{equation}
\left | \beta_{4K - \ell}(K) \right |
\leq \left ( 4K \right )^\ell 
\cdot
\left ( 2K-1 \right )^{\ell}
\leq
\left ( 8K^2 \right )^\ell .
\label{eq:Bound+Z}
\end{equation}

Collecting the bounds (\ref{eq:Bound+X}), (\ref{eq:Bound+Y})
and (\ref{eq:Bound+Z}) we see from 
(\ref{eq:InequalityForAbsolutevalue}) that
\begin{eqnarray}
\left | a_{\ell} (K) \right |
&\leq&
\sum_{k=0}^4
\frac{K^{2k}}{k!} \cdot \left ( 12 K^2 \right )^{\ell - k}
+ \left ( 8K^2 \right )^\ell .
\nonumber \\
&=& C_\ell \left ( 12 K^2 \right )^\ell
\end{eqnarray}
with
\[
C_\ell 
:= 
\sum_{k=0}^4 \frac{1}{ k! \cdot 12^k } 
+ 
\left ( \frac{8}{12} \right )^\ell .
\]
It is a simple matter to check that $C_\ell \leq 2$.
\myendpf

\subsection{A proof of Lemma \ref{lem:first_six_coef}}

The basis for the proof can be found in the expression
(\ref{eq:CoefficientFormula+Equiv}) for the coefficients
$a_{0} (K) , \ldots , a_{4K} (K)$.

For $\ell=0$, this expression becomes
\[
a_{0} (K)
= \beta_{0,4K}(K) - \beta_{4K}(K) 
= 1 - 1 = 0.
\]

For $\ell=1$, we get
\begin{eqnarray}
a_1(K)
&=& \beta_{0,4K-1}(K) + K^2 \beta_{1,4K-1} (K)
- \beta_{4K-1} (K)
\nonumber \\
&=& 
- \sum_{\ell=1}^{K} (3K-\ell) 
- \sum_{j=1}^{3K} (3K -j) 
+ K^2 
- \left ( - 4 \sum_{i=K}^{2K-1} i \right ) = 0
\nonumber
\end{eqnarray}
where we have used the formula
(\ref{eq:FromRootsToCoefficients}) to evaluate 
$\beta_{0,4K-1}(K)$ and $\beta_{4K-1} (K)$.

For $\ell=2$, this approach yields
\begin{eqnarray}
\lefteqn{a_2(K)} & &
\nonumber \\
&=& \beta_{0,4K-2}(K) + K^2 \beta_{1,4K-2} (K)
+  \frac{K^2(K-1)^2}{2} \beta_{2,4K-2} (K)
- \beta_{4K-2}(K) 
\nonumber \\
&=& \beta_{0,4K-2}(K) + K^2 \beta_{1,4K-2} (K)
+  \frac{K^2(K-1)^2}{2}
- \beta_{4K-2}(K) 
\nonumber
\end{eqnarray}
and this leads to
\begin{eqnarray}
a_2(K)
&=&
\sum_{i=1}^{3K-2}i\sum_{j=i+1}^{3K-1}j 
+ \sum_{i=2K}^{3K-2}i\sum_{j=i+1}^{3K-1}j 
+ \left(\sum_{i=2K}^{3K-1}i\right)\left(\sum_{i=1}^{3K-1}i\right)
\\\nonumber
& & ~ -K^2\left(\sum_{i=1}^{3K-2}i+\sum_{i=2K-1}^{3K-2}i\right)
+
\frac{K^2 (K-1)^2}{2}
\\\nonumber
& & ~ -\left( {4 \choose 1 } ^ {2} \cdot
\sum_{i=K}^{2K-2}i\sum_{j=i+1}^{2K-1}j 
+ { 4 \choose 2 } \sum_{i=K}^{2K-1}i^2\right)
\\\nonumber
&=& 0.
\end{eqnarray}

For $\ell=3$, straightforward computations give
\begin{eqnarray}
\lefteqn{a_3(K)}&&
\\\nonumber
&=&-\left(\sum_{v=1}^{3K -3}v
\sum_{i=v+1}^{3K-2}i\sum_{j=i+1}^{3K-1}j + \sum_{v=2K}^{3K -3} v
\sum_{i=v+1} ^ {3K-2}i\sum_{j=i+1}^{3K-1}j \right)
\\\nonumber
& & ~ -\left(\sum_{i=1}^{3K -1} i \cdot \sum_{i=2K} ^ {3K-2} i
\sum_{j=i+1} ^ {3K-1} j + \sum_{i=2K} ^ {3K-1} i \cdot  \sum_{i=1}
^ {3K-2} i \sum_{j=i+1} ^ {3K-1} j \right)
\\\nonumber
& & ~ + K^2\left( \sum_{i=1}^{3K-3}i\sum_{j=i+1}^{3K-2}j+
\sum_{i=2K-1}^{3K-3}i\sum_{j=i+1}^{3K-2}j +
\left(\sum_{i=2K-1}^{3K-2}i\right)\left(\sum_{i=1}^{3K-2}i\right)
\right)
\\\nonumber
& & ~ - \frac{K^2
(K-1)^2}{2}\left(\sum_{i=1}^{3K-3}i+\sum_{i=2K-2}^{3K-3}i\right)+\frac{K^2
(K-1)^2(K-2)^2}{6}
\\\nonumber
& & ~ + {4\choose 1} ^3 \cdot \sum_{v=K}^{2K -3}v
\sum_{i=v+1}^{2K-2}i\sum_{j=i+1} ^ {2K-1} j + {4 \choose 2}{4
\choose 1} \sum_{j=K}^{2K-1} j^2 \left( \sum_{i=K} ^ {2K-1} i - j
\right)
\\ \nonumber
& & ~ + { 4 \choose 3 } \sum_{i=K}^{2K-1}i^3
\\\nonumber
&=&K^4 
\end{eqnarray}
as announced.

For $a_4(K)$, we
proceed in a similar manner to get\footnote{Evaluating the
expression (\ref{eq:a_4_eval}) (as well as (\ref{eq:a_5_eval})
given next) by hand is quite cumbersome. To avoid this, one can
make use of a computer software (e.g., Mathematica, MATLAB)
that can perform computations symbolically.}
\begin{eqnarray}\label{eq:a_4_eval}
\lefteqn{a_4(K)}&&
\\\nonumber
&=& \sum_{l=1} ^ {3K-4}l \sum_{v = l+1} ^ {3K -3} v
\sum_{i=v+1}^{3K-2}i\sum_{j=i+1}^{3K-1}j + \sum_{l=2K}^{3K -4} l
\sum_{v=l+1} ^ {3K -3}v \sum_{i=v+1}^{3K-2}i\sum_{j=i+1}^{3K-1}j
\\\nonumber
& & ~ + \sum_{i=1} ^ {3K -1} i \cdot \sum_{v=2K} ^ {3K-3} v
\sum_{i=v+1} ^ {3K-2} i \sum_{j=i+1} ^ {3K-1} j + \sum_{i=2K}^
{3K-1} i \cdot \sum_{v=1}^{3K-3} v \sum_{i=v+1} ^ {3K-2} i
\sum_{j=i+1} ^ {3K-1} j
\\\nonumber
& & ~ + \left ( \sum_{i=1}^{3K-2} i \sum_{j=i+1} ^ {3K-1} j
\right) \left ( \sum_{i=2K}^{3K-2} i \sum_{j=i+1} ^ {3K-1} j
\right) - K^2 \sum_{v=1}^{3K-4}v \sum_{i=v+1} ^ {3K-3} i
\sum_{j=i+1}^{3K-2}j
\\\nonumber
& & ~ - K^2 \left( \sum_{v=2K-1} ^ {3K-4} v \sum_{i=v+1} ^ {3K-3}
i \sum_{j=i+1}^{3K-2}j + \sum_{i=1} ^ {3K-2} i \cdot
\sum_{i=2K-1}^{3K-3}i \sum_{j=i+1} ^ {3K-2} j \right )
\\\nonumber
& & ~ -K^2 \sum_{i=2K-1} ^ {3K-2} i \cdot \sum_{i=1}^{3K-3}i
\sum_{j=i+1} ^ {3K-2} j  + \frac{K^2 (K-1)^2}{2}
\sum_{i=1}^{3K-4}i \sum_{j=i+1} ^ {3K-3} j
\\\nonumber
& & ~ + \frac{K^2 (K-1)^2}{2}\left( \sum_{i=2K-2}^{3K-4}i
\sum_{j=i+1} ^ {3K-3} j + \left( \sum_{i=2K-2} ^ {3K-3} i\right)
\left( \sum_{j=1} ^ {3K-3} j \right) \right)
\\\nonumber
& & ~ - \frac{K^2 (K-1)^2(K-2)^2}{6}\left(\sum_{i=1} ^ {3K-4}i +
\sum_{i=2K-3} ^ {3K-4} i \right)
\\\nonumber
& & ~ + \frac{K^2 (K-1)^2(K-2)^2 (K-3)^2}{24} - {4\choose1}^4
\cdot \sum_{l=K} ^ {2K -4} l \sum_{v=K+1} ^ {2K -3} v \sum_{i=v+1}
^ {2K-2} i \sum_{j=i+1} ^ {2K-1} j
\\\nonumber
& & ~ - {4 \choose 2}{4 \choose 1}^2 \cdot \sum_{v=K} ^ {2K-1} v^2
\left(\sum_{i=K} ^ {2K-2} i \sum_{j=i+1} ^ {2K-1} j - v \sum_{i=K}
^ {2K-1} i + v^2 \right)
\\\nonumber
& & ~ - {4 \choose 3}{4 \choose 1} \sum_{j=K} ^ {2K-1} j^3
\left(\sum_{i=K} ^ {2K-1}i - j \right) -  { 4 \choose 2 } ^ 2
\cdot \sum_{i=K} ^ {2K-2} i ^ 2 \sum_{ j= i + 1 } ^ { 2 K - 1} j ^
2 - \sum_{ i = K } ^ { 2 K - 1 } i ^ 4
\\\nonumber
&=&-6K^6+6K^5-K^4.
\end{eqnarray}

Finally, $a_5(K)$ is given by
\begin{eqnarray}\label{eq:a_5_eval}
\lefteqn{a_5(K)}&&
\\\nonumber
&=& -\sum_{u=1} ^ {3K-5} u \sum_{l=u+1} ^ {3K-4} l \sum_{v = l+1}
^ {3K -3} v \sum_{i=v+1}^{3K-2}i\sum_{j=i+1}^{3K-1}j
\\\nonumber
& & ~ - \sum_{u=2K}^{3K -5} u \sum_{l=u+1}^{3K -4} l \sum_{v=l+1}
^ {3K -3}v \sum_{i=v+1}^{3K-2}i\sum_{j=i+1}^{3K-1}j
\\\nonumber
& & ~ - \sum_{i=1} ^ {3K -1} i \cdot \sum_{l=2K} ^ {3K-4} l
\sum_{v=l+1} ^ {3K-3} v \sum_{i=v+1} ^ {3K-2} i \sum_{j=i+1} ^
{3K-1} j
\\ \nonumber
& & ~ -  \sum_{i=2K}^ {3K-1} i \cdot \sum_{l=1}^{3K-4} l
\sum_{v=l+1} ^ {3K-3} v \sum_{i=v+1} ^ {3K-2} i \sum_{j=i+1} ^
{3K-1} j
\\\nonumber
& & ~ - \left ( \sum_{v=1}^{3K-3} v \sum_{i=v+1} ^ {3K-2} i
\sum_{j=i+1} ^ {3K-1} j \right) \left ( \sum_{i=2K}^{3K-2} i
\sum_{j=i+1} ^ {3K-1} j \right)
\\\nonumber
& & ~ - \left ( \sum_{v=2K}^{3K-3} v \sum_{i=v+1} ^ {3K-2} i
\sum_{j=i+1} ^ {3K-1} j \right) \left ( \sum_{i=1}^{3K-2} i
\sum_{j=i+1} ^ {3K-1} j \right)
\\ \nonumber
& & ~ + K^2 \left(\sum_{l=1} ^ {3K-5} l \sum_{v=l+1} ^ {3K-4} v
\sum_{i=v+1} ^ {3K-3} i \sum_{j=i+1}^{3K-2}j
 + \sum_{l=2K-1} ^ {3K-5} l \sum_{v=l+1} ^ {3K-4} v
\sum_{i=v+1} ^ {3K-3} i \sum_{j=i+1}^{3K-2}j \right)
\\\nonumber
& & ~ + K^2 \sum_{i=1} ^ {3K-2} i \cdot \sum_{v=2K-1}^{3K-4} v
\sum_{i=v+1} ^ {3K-3} i \sum_{j=i+1} ^ {3K-2} j
\\\nonumber
& & ~ + K^2 \sum_{i=2K-1} ^ {3K-2} i \cdot \sum_{v=1} ^ {3K-4} v
\sum_{i=v+1} ^ {3K-3} i \sum_{j=i+1} ^ {3K-2} j
\\\nonumber
& & ~ + K^2 \left ( \sum_{i=1}^{3K-3} i \sum_{j=i+1} ^ {3K-2} j
\right) \left ( \sum_{i=2K-1} ^ {3K-3} i \sum_{j=i+1} ^ {3K-2} j
\right)
\\\nonumber
& & ~ - \frac{K^2 (K-1)^2}{2} \left( \sum_{v=1} ^ {3K-5} v
\sum_{i=v+1} ^ {3K-4} i \sum_{j=i+1} ^ {3K-3} j + \sum_{v=2K-2} ^
{3K-5} v \sum_{i=v+1} ^ {3K-4} i \sum_{j=i+1} ^ {3K-3} j \right)
\\\nonumber
& & ~ - \frac{K^2 (K-1)^2}{2}\left( \sum_{i=1} ^ {3K-3} i \cdot
\sum_{i=2K-2}^{3K-4}i \sum_{j=i+1} ^ {3K-3} j + \sum_{i=2K-2} ^
{3K-3} i \cdot \sum_{i=1}^{3K-4}i \sum_{j=i+1} ^ {3K-3} j \right)
\\\nonumber
& & ~ + \frac{K^2 (K-1)^2(K-2)^2}{6}\left(\sum_{i=1} ^ {3K-5} i
\sum_{j=i+1} ^ {3K-4} j + \sum_{i=2K-3} ^ {3K-5} i \sum_{j=i+1} ^
{3K-4} j \right)
\\\nonumber
& & ~ +  \frac{K^2 (K-1)^2(K-2)^2}{6} \sum_{i=1} ^ {3K-4} i \cdot
\sum_{i=2K-3} ^ {3K-4} j
\\ \nonumber
& & ~ - \frac{K^2 (K-1)^2(K-2)^2 (K-3)^2}{24}\left( \sum_{i=1} ^
{3K-5} i + \sum_{i=2K-4} ^ {3K-5} i \right)
\\\nonumber
& & ~ + {4\choose1}^5 \cdot \sum_{u=K} ^ {2K -5} u \sum_{l=u+1} ^
{2K -4} l \sum_{v=K+1} ^ {2K -3} v \sum_{i=v+1} ^ {2K-2} i
\sum_{j=i+1} ^ {2K-1} j
\\\nonumber
& & ~ + {4 \choose 2} {4 \choose 1}^3
\\\nonumber
& & ~ ~ ~ \times \sum_{l=K} ^{2K-1} l^2 \left( \sum_{v=K} ^{2K-3}
v \sum_{i=m+1} ^ {2K-2} i \sum_{j=i+1} ^ {2K-1} j - l \sum_{i=K} ^
{2K-2} i \sum_{j=i+1} ^ {2K-1} j + l^2 \sum_{i=K} ^ {2K-1} -l^3
\right)
\\\nonumber
& & ~ + {4 \choose 2}^2{4 \choose 1} \cdot \sum_{v=K} ^ {2K-1} v
\left(\sum_{i=K} ^ {2K-2} i^2 \sum_{j=i+1} ^ {2K-1} j^2 - v^2
\sum_{i=K} ^ {2K-1} i^2 + v^4 \right)
\\\nonumber
& & ~ + {4 \choose 3}{4 \choose 2} \sum_{i=K} ^ {2K-1} i^3
\left(\sum_{j=K} ^ {2K-1}j^2 - i^2 \right)
\\\nonumber
& & ~ +  {4 \choose 3} {4 \choose 1}^2 \cdot \sum_{v=K}^{2K-1} v^3
\left( \sum_{i=K} ^ {2K-2} i \sum_{j=i+1} ^ {2K-1} j- v \sum_{i=K}
^ {2K-1} +v^2 \right)
\\\nonumber
& & ~ +  {4 \choose 4} {4 \choose 1} \sum_{i=K} ^{2K-1} i^4 \left(
\sum_{j=K} ^ {2K-1} j - i \right)
\\\nonumber
&=&-\frac{1}{120}K^{10}+\frac{1}{6}K^9+\frac{199}{12}K^8-34K^7+\frac{1207}{120}K^6+\frac{161}{6}K^5
\\\nonumber
& & ~ - \frac{209}{6}K^4+20K^3-\frac{24}{5}K^2.
\end{eqnarray}
\myendpf

\end{document}